\documentclass[10pt,a4paper,reqno]{amsart}
\usepackage{amsfonts,amsthm,latexsym,amsmath,amssymb,amscd,epsfig}
\usepackage{graphics,graphicx}
\newtheorem{Theorem}{Theorem}
\newtheorem{Lemma}{Lemma}
\newtheorem{Cor}{Corollary}
\newtheorem{Prop}{Proposition}

\newtheorem{Claim}{Claim}
\newtheorem{theo+}           {Theorem}
\newtheorem{prop+}           {Proposition}
\newtheorem{coro+}           {Corollary}
\newtheorem{lemm+}           {Lemma}
\newtheorem{claim+}          {Claim}

\theoremstyle{definition}
\newtheorem{defi+}           {Definition}

\newtheorem*{m-p}            {Main Problem}
\newtheorem{not+}            {Notation} 
\newtheorem{Ex}              {Example}
\newtheorem{Prob}            {Problem}

\theoremstyle{remark}
\newtheorem{rema+}           {Remark}

\newenvironment{remark}{\begin{rema+}}{\end{rema+}}
\newenvironment{definition}{\begin{defi+}}{\end{defi+}}
\newenvironment{notation}{\begin{not+}}{\end{not+}}

\newcommand{\al}{\alpha}
\newcommand{\be}{\beta}

\newcommand{\te}{\theta}
\newcommand {\bC} {\mathbb C}
\newcommand {\bR} {\mathbb R}

\newcommand{\bN}{\mathbb N}
\newcommand {\D} {\mathcal D}

\newcommand {\Ga} {\Gamma}
\newcommand {\De}{\Delta}
\newcommand {\La}{\Lambda}
\newcommand {\Om}{\Omega}
\newcommand {\ga} {\gamma}
\newcommand {\eps} {\epsilon}
\newcommand {\de} {\delta}
\newcommand{\ze}{\zeta}
\newcommand {\na}{\nabla}
\newcommand{\vf}{\varphi}
\newcommand{\cD}{\mathcal{D}}
\newcommand{\cM}{\mathcal{M}}
\newcommand{\cA}{\mathcal{A}}
\newcommand{\Si}{\Sigma}
\newcommand{\pa}{\partial}

\begin{document}
\numberwithin{equation}{section}

\title[Piecewise harmonic functions]{Piecewise harmonic subharmonic functions\\
and positive Cauchy transforms}
\author[J.~Borcea]{Julius Borcea}
\address{Department of Mathematics, Stockholm University, SE-106 91, Stockholm,
Sweden}
\email{julius@math.su.se}
\author[R.~B\o gvad]{Rikard B\o gvad}
\address{Department of Mathematics, Stockholm University, SE-106 91, Stockholm,
Sweden}
\email{rikard@math.su.se}

\begin{abstract}
We give a local characterization of the class of functions having positive 
distributional derivative with respect to $\bar{z}$ that are almost everywhere
equal to one of finitely many analytic functions and satisfy some mild 
non-degeneracy assumptions. As a consequence, we give conditions that 
guarantee that any subharmonic piecewise harmonic function 
coincides locally with the 
maximum of finitely many harmonic functions and we describe 
the topology of their level curves. These results are valid in a quite general 
setting as they assume no {\em \`a priori} conditions on the differentiable 
structure of the support of the associated Riesz measures. We also discuss 
applications to positive Cauchy transforms and we consider several examples 
and related problems.
\end{abstract}

\maketitle

\tableofcontents

\section{Introduction}\label{s-1}

\enlargethispage{0.2cm}

One of the most frequently used constructions in complex analysis and geometry
is to consider the maximum of a finite number of pairwise distinct harmonic 
functions. As is well known, the result is a subharmonic function which is 
also piecewise harmonic. A quite natural problem is to investigate the 
converse direction, namely study the class of functions generated by this basic
albeit fundamental procedure. Its classical flavor \cite{Ha} and some 
important applications -- some of which are listed below -- further motivate
a deeper study of this question on which surprisingly little seems to be known.
 In this paper we answer this question by giving a local 
characterization of the aforementioned class of functions in generic cases and 
in the process we establish several remarkable properties for this class. 
In particular, we show that any subharmonic piecewise harmonic function may 
essentially be realized as the maximum of finitely many harmonic functions. 

\subsection{Piecewise Harmonic and Piecewise Analytic Functions}\label{s-11}

Let us first define a fairly general notion.

\begin{definition}\label{d-r1}
Let $X$ be a real or complex subspace of the space of smooth functions in
a domain (open connected set) $U$ in $\bR^2$ or $\bC$. We say that a function
$\varphi$ is {\em piecewise in} $X$ if one can find finitely many
pairwise disjoint open sets $M_i$, $1\le i\le r$, in $U$ and pairwise distinct
functions $\varphi_i\in X$, $1\le i\le r$, such that
\begin{itemize}
\item[(i)] $\varphi=\varphi_i$ in $M_i$, $1\le i\le r$;
\item[(ii)] $U\setminus \bigcup_{i=1}^{r}M_i$ is of Lebesgue measure $0$.
\end{itemize}
The set of all functions that are piecewise in $X$ is denoted by $PX$.
\end{definition}

\begin{remark}\label{r-ex-px}
It is not difficult to see that $PX$ is actually a (real or complex) vector
space. This as well as further properties of $PX$ functions and related 
concepts are discussed in the Appendix. 
\end{remark}

Note that since $PX$ functions are locally integrable they define
distributions and their derivatives are therefore defined in the distribution
sense (and functions are identified if they define the same distributions). 
In particular, if $\varphi\in PX$ one can form $\De \varphi\in\cD'(U)$ and 
also $\partial_z\varphi,\partial_{\bar{z}}\varphi\in\cD'(U)$ if $X$ is complex.

We now specialize Definition~\ref{d-r1} to obtain the main objects of
our study, namely the spaces of {\em piecewise harmonic} and
{\em piecewise analytic functions}, respectively.

\begin{notation}\label{not-r1}
Fix a domain $U\subset \bC$, let $H=H(U)$ be the real space of
(real-valued) harmonic functions in $U$ and $A=A(U)$ be the complex 
space of analytic functions in $U$. By Definition~\ref{d-r1} the following
holds:
\begin{itemize}
\item[(a)] Given a 
piecewise harmonic function $\varphi\in PH$ there exists a finite family 
of pairwise disjoint open sets $\{M_i\}_{i=1}^{r}$ in $U$ covering $U$ up to 
a set of Lebesgue measure $0$ and a corresponding family 
of pairwise distinct harmonic functions $\{H_{i}(z)\}_{i=1}^{r}$ in $U$ such 
that
\begin{equation}\label{eq:MDH}
\varphi(z)=\sum_{i=1}^{r}H_{i}(z)\chi_{i}(z)
\end{equation} 
a.e.~in $U$, where $\chi_{i}$ is the characteristic function of the set 
$M_{i}$, $1\le i\le r$; 
\item[(b)] Similarly, any piecewise analytic function $\Phi\in PA$
may be represented as
\begin{equation}\label{eq:MDA}
\Phi(z) =\sum_{i=1}^{r}A_{i}(z)\chi_{i}(z)
\end{equation}
a.e.~in $U$, where $M_i$ and $\chi_i$, $1\le i\le r$, are as in (a) and 
$\{A_{i}(z)\}_{i=1}^{r}$ is a family of pairwise distinct analytic functions 
in $U$. Given this data and a point $p\in U$ we set
\begin{equation}\label{a-h}
H_i(z)=\Re\left[\int_{p}^z A_{i}(w)dw\right],\quad z\in U,\,1\le i\le r.
\end{equation}
These are well-defined harmonic functions in $U$ provided that $U$ is simply 
connected, which we tacitly assume throughout unless otherwise stated.
\end{itemize}
\end{notation}

We stress the fact that in the above definitions no regularity ($C^1$) 
conditions are assumed on the negligible set 
$U\setminus \bigcup_{i=1}^{r}M_i$. 
Note also that Definition \ref{d-r1} and Notation \ref{not-r1} are 
merely a convenient way of saying 
that a $PH$ function $\phi$ equals one of finitely many harmonic functions in 
certain prescribed sets. 
Therefore $PH$ functions need not be continuous nor subharmonic and one can 
hardly expect any interesting statements in this kind of generality.
The same philosophy applies to $PA$ functions: as defined above,
a function $\Phi$ is $PA$ if it is equal to one of finitely many analytic 
functions in certain open sets. Thus $PA$ functions need not be continuous 
and this will not be case either in our situation.

%Let us first make the following definition.
%\begin{definition}\label{def-ph}
%A function $\vf$ defined in a domain $U\subset \bC$ is called 
%{\em piecewise harmonic}, or $PH$ for short, if there exists a finite family 
%of  and 
%an associated family of disjoint open sets $\{M_i\}_{i=1}^{r}$ with 
%$M_i\subset $,
%$1\le i\le r$, such that $\{M_i\}_{i=1}^{r}$ 
%forms a covering of $U$ up to a set of Lebesgue measure $0$ and
%$$\varphi(z)=\sum_{i=1}^{r}H_{i}(z)\chi_{i}(z)\text{ in }U,$$
%where $\chi_{i}$ is the characteristic function of the set 
%$M_{i}$, $1\le i\le r$.
%\end{definition}

\subsection{Canonical Piecewise Decompositions}\label{s-12}

Note that conditions (i)--(ii) in Definition \ref{d-r1} remain valid if 
non-empty Lebesgue negligible sets are subtracted from the sets $M_i$, 
so it is in general impossible to say something about the boundaries 
of these sets. However, the inclusions 
$M_i\subseteq U\setminus\text{supp}(\varphi-\varphi_i)$, $1\le i\le r$, 
always hold, where the 
supports are defined in the distribution sense (recall from \S \ref{s-11} that 
$PX$ functions are locally integrable and $L_{loc}^1(U)$ is viewed as a 
subspace of $\cD'(U)$). 
Now both $X=H(U)$ and $X=A(U)$ are examples of function 
spaces satisfying the {\em unique continuation property}, i.e., $f\equiv 0$ 
in $U$ if $f\in X$ vanishes in some open non-empty subset of $U$. 
In view of the above inclusions, for spaces with this property 
one can reformulate Definition \ref{d-r1} in a more canonical 
way as follows. 

\begin{definition}\label{d-r2}
Let $X$ be a real or complex subspace of the space of smooth functions in
a domain $U$ in $\bR^2$ or $\bC$. Assume that $X$ satisfies the unique 
continuation property and let $\varphi\in L_{loc}^1(U)$. Then $\varphi\in PX$ 
($\varphi$ is piecewise in $X$) if one can find pairwise distinct 
elements $\varphi_i\in X$, $1\le i\le r$, such that the set 
$\Gamma:=\bigcap_{1\le i\le r}\text{supp}(\varphi-\varphi_i)$ is of Lebesgue 
measure $0$.
\end{definition}

Setting $M_i=U\setminus\text{supp}(\varphi-\varphi_i)$, $1\le i\le r$, in 
Definition \ref{d-r2} we see that $M_i$ is the largest open 
set in which $\varphi-\varphi_i$ vanishes (as a distribution or almost 
everywhere). Further useful properties of the canonical piecewise 
decomposition 
of the $PX$ function $\varphi$ given in Definition \ref{d-r2} are gathered 
in the next lemma. Henceforth by a ``continuous function'' we mean a function 
in $L^1_{loc}(U)$ which agrees almost everywhere with a continuous function 
in $U$.

\begin{Lemma}\label{l-can}
In the above notation the following holds:
\begin{itemize}
\item[(i)] $\bigcup_{1\le i\le r}M_i=U\setminus \Gamma$; 
\item[(ii)] $\overline{M}_i\cap M_j=\emptyset$, $1\le i\neq j\le r$; 
\item[(iii)] $M_i=\mathring{\overline{M}_i}$ (i.e., $M_i$ is the interior 
of $\overline{M}_i$), $1\le i\le r$; 
\item[(iv)] $\Gamma=\bigcup_{1\le i<j\le r}\overline{M}_i\cap \overline{M}_j$.
\item[(v)] If $\varphi$ is continuous then $\Gamma\subseteq g^{-1}(0)$, where 
%is contained in the zero set of the function 
$g:=\prod_{1\le i<j\le r}(\varphi_i-\varphi_j)$.
\end{itemize}
\end{Lemma}

\begin{proof}
The first statement is obviously true by the (canonical) 
definition of the sets $M_i$, $1\le i\le r$. To prove (ii) suppose that 
$i\neq j$ and 
$p\in \overline{M}_i\cap M_j$. Then one can find $q\in M_i$ arbitrarily close 
to $p$ with $q\in M_i\cap M_j$. Since $q\notin \text{supp}(\varphi-\varphi_i)$
and $q\notin \text{supp}(\varphi-\varphi_j)$ one gets 
$q\notin \text{supp}(\varphi_i-\varphi_j)$ and the unique continuation 
property implies that $\varphi_i=\varphi_j$, which contradicts the fact that 
$\varphi_i\neq \varphi_j$.

To show (iii) assume that $p\in\mathring{\overline{M}_i}$. Then 
there exists an (open) neighborhood $N$ of $p$ which is contained in 
$\overline{M}_i$. Since $\overline{M}_i\cap M_j=\emptyset$ if $j\neq i$ 
(cf.~(ii)) it follows that $N\subset M_i\cup \Gamma$. 
Hence $\varphi=\varphi_i$ 
in $N$ and $N\subset M_i$, so that in particular $p\in M_i$.

Clearly, $\bigcup_{1\le i\le r}\overline{M}_i=U$. 
Therefore, if $p\in \Gamma$ then 
$p\in \overline{M}_i$ for some $i$ and $p$ must then be a boundary point of 
$M_i$. Assume that $p\notin \overline{M}_j$ whenever $j\neq i$. Then there is 
a neighborhood $N$ of $p$ such that $N\cap M_j=\emptyset$ for $j\neq i$. Hence 
$N\subset \overline{M}_i$ and it follows from (iii) that 
$p\in\mathring{M_i}$. This gives a contradiction (since $p$ is a boundary 
point 
of $M_i$) and shows that $p\in \overline{M}_i\cap \overline{M}_j$ for 
some $j\neq i$, which proves (iv).

Finally, if $\varphi$ is continuous then 
$\varphi=\varphi_i$ in $\overline{M}_i$ and 
$\varphi=\varphi_j$ in $\overline{M}_j$ hence $\varphi_i=\varphi_j$ 
in $\overline{M}_i\cap \overline{M}_j$ and thus $g\equiv 0$ in 
$\overline{M}_i\cap \overline{M}_j$ for $i\neq j$, so that by (iv) 
$g\equiv 0$ in $\Gamma$.
\end{proof}

The familiar ``maximum construction'' that we alluded to at the beginning of 
this introduction yields natural examples of $PH$ and $PA$ functions. We 
recall briefly the interplay between the classes of functions obtained in 
this case: 

\begin{Ex}\label{ex-bla}
Let $\{H_{i}(z)\}_{i=1}^{r}$ be a finite family of pairwise distinct harmonic 
functions in a domain $U\subset \bC$. Then 
$\varphi(z):=\max_{1\le i\le r} H_{i}(z)$ is a (subharmonic) 
$PH$ function. Indeed, set
$\Om:=\{z\in U\mid H_{k}(z)\neq H_{l}(z), 1\le k\neq l\le r\}$, 
let $M_{i}$ be the (open) set consisting of 
those $z\in \Om$ for which $\varphi(z)=H_{i}(z)$ and denote by $\chi_{i}$ 
the characteristic function of $M_{i}$, $1\le i\le r$. It is clear that 
$U\setminus \Om$ is Lebesgue negligible, so that $\{M_i\}_{i=1}^{r}$ 
forms a covering of $U$ up to a set of Lebesgue measure $0$ and
\begin{equation*}
\varphi(z)=\sum_{i=1}^{r}H_{i}(z)\chi_{i}(z)
\end{equation*} 
a.e.~in $U$. Moreover, the subharmonicity of $\varphi$ implies that 
$\nu:=\partial^2 \varphi/\partial \bar z\partial z\ge 0$  
in the sense of distributions. In fact $\nu$ 
is a positive measure supported on the (finite) union of level curves 
$\{z\in U\mid H_{i}(z)-H_{j}(z)=0\}$, $1\le i\neq j\le r$. One can show 
that in this case the support actually determines the measure 
(Theorem \ref{coro6} in \S \ref{s-2}). 

Now the derivative of $\vf$, again in the distribution sense, 
inherits a similar property only this time with respect to analytic 
functions. Classical results yield namely 
\begin{equation*}
    \partial \varphi(z)/\partial z=
\sum_{i=1}^{r}A_{i}(z)\chi_{i}(z)
\end{equation*}
a.e.~in $U$, where $A_{i}:=\partial H_{i}/\partial z$, $1\le i\le r$, 
are analytic 
functions in $U$ (cf.~Proposition~\ref{coro5} in \S \ref{s-2}). 
Hence $\partial \varphi(z)/\partial z$ is a $PA$ function.
Note that the above relation may be 
reformulated as saying that $\varphi$ satisfies a.e.~in $U$ 
the differential equation
$P(\partial \varphi(z)/\partial z,z)=0$, 
where $P(y,z):=\prod_{i=1}^{r}(y-A_{i}(z))$ is a polynomial in $y$ with 
coefficients that are holomorphic in $U$. 
\end{Ex}

%\begin{definition}\label{def-pa}
%A function $\Phi$ defined in a domain $U\subset \bC$ is called 
%{\em piecewise analytic}, or $PA$ for short, if there exists a finite family 
%of pairwise distinct analytic functions $\{A_{i}(z)\}_{i=1}^{r}$ in $U$ and 
%an associated family of disjoint open sets $\{M_i\}_{i=1}^{r}$ with
%$M_i\subset \Om:=\{z\in U\mid A_{k}(z)\neq A_{l}(z), 1\le k\neq l\le r\}$, 
%$1\le i\le r$, such that $\{M_i\}_{i=1}^{r}$ forms a covering of $U$ up to a 
%set of Lebesgue measure $0$ and 
%$$\Phi(z)=\sum_{i=1}^{r}A_{i}(z)\chi_{i}(z)\text{ in }U,$$
%where $\chi_{i}$ is the characteristic function of the set 
%$M_{i}$, $1\le i\le r$.
%\end{definition}

\subsection{Main Problem and Results}\label{s-13}

$PA$ functions occur naturally -- and this was our original 
motivation -- in various contexts, such as the study of the 
asymptotic behavior of polynomial solutions to ordinary differential 
equations \cite{BR,BBS,Fe,Wa1}, the theory of Stokes lines \cite{K,Wa2} and
orthogonal polynomials \cite{DZ}. 
In the aforementioned contexts $PA$ 
functions are mostly constructed as limits and thus one has no control 
on the differentiable structure of the resulting sets $M_{i}$.  
It is therefore important to describe the local and global 
structure of $PA$ functions both with and without additional regularity 
assumptions -- such as piecewise $C^1$-boundary conditions on the 
sets $M_{i}$, see \S \ref{s-2} -- and this is the primary objective of this
paper. To state our main problem it is convenient to use:

\begin{notation}\label{not-Sigma}
Given a domain $U\subset \bC$ let 
$\Si(U)=\{f\in\cD'(U)\mid \partial_{\bar{z}}f\ge 0\}$.
\end{notation}

Clearly, $\pa_z\varphi\in\Si(U)$ if $\varphi$ is subharmonic in $U$,
which holds e.g.~for the maximum of finitely many harmonic 
functions. For a (known) converse see the Appendix.

\begin{m-p}
Let $\Phi\in\Si(U)$ be a $PA$ function in a given domain $U\subset \bC$. Find 
conditions that guarantee that $\Phi$ is locally (or globally) of the form 
$\pa_{z}\varphi$, where $\varphi$ is the maximum of a finite number of 
harmonic functions in $U$.
\end{m-p}

The necessity of assuming $\partial_{\bar{z}}\Phi\ge 0$ in the 
Main Problem will soon become quite clear and is further illustrated in 
Example \ref{ex-bla2}, see also Lemma \ref{ref-l-bla} in the Appendix.
We give four answers to the above problem which may be summarized 
(in terms of the mutual implications among them) as follows:
\begin{equation}\label{eq-chain}
\text{Theorem \ref{main-r1}}\Longrightarrow  
\text{Corollary \ref{Cor5}}\Longrightarrow
\text{Corollary \ref{cor-ip-conv}}\Longrightarrow
\text{Theorem \ref{th-ref}}.
\end{equation}

We formulate here just the first 
(Theorem \ref{main-r1}) and third 
(Corollary \ref{cor-ip-conv}) main results of this paper. The fourth one 
(Theorem \ref{th-ref}) is an alternative approach to the Main Problem 
suggested by our referee, as were several ideas used in this paper.

\begin{Theorem}\label{main-r1} 
Let $\Phi\in\Si(U)$ be a $PA$ function as in \eqref{eq:MDA} and assume 
that $p\in U$ satisfies the following conditions:
\begin{itemize}
\item[(i)] $p\in\overline{M}_i$, $1\le i\le r$;
%$N(p)\cap M_{i}$ has positive Lebesgue measure for any 
%neighborhood $N(p)$ of $p$ and $1\le i\le r$; 
\item[(ii)] $A_{i}(p)-A_{k}(p)\notin \bR(A_{j}(p)-A_{k}(p))$ for any triple
of distinct indices $(i,j,k)$ in $\{1,\ldots,r\}$.
\item[(iii)] $A_{i}(p)\neq A_{k}(p)$ for any pair of distinct indices
$(i,k)$ in $\{1,\ldots,r\}$.
\end{itemize}
There exists a
neighborhood $\widetilde N(p)$ of $p$ such that
$\Phi=2\partial \varphi /\partial z$ a.e.~in $\widetilde N(p)$, 
where 
$\varphi(z)=\max_{1\le i\le r}H_i(z)$ and the $H_i$'s are the harmonic 
functions defined in~\eqref{a-h}.
\end{Theorem}

A word about each of the three 
conditions imposed in Theorem \ref{main-r1} is in 
order:
\begin{itemize}
\item[(i)] suggests defining the following index set for any $p\in U$:
\begin{equation}\label{ip}
I(p)=\{j\in\{1,\ldots,r\}\mid p\in\overline{M}_j\}
\end{equation}
and $i(p)=|I(p)|$. Condition (i) then requires that $i(p)=r$, i.e., every 
set $M_i$ is ``active''. This will be tacitly assumed throughout;
\item[(ii)] is the most important assumption and amounts to the 
requirement that for all distinct indices $i,j,k\in\{1,\ldots,r\}$ the level 
curves $H_i=H_k$ and $H_j=H_k$ should meet transversally at $p$ (i.e., 
the critical sets $\Gamma_{i,j,k}$ defined in
\eqref{eq-crit-set} below do not contain $p$). For an illustration of the 
necessity of this assumption 
see Example~\ref{ex1} and Figure 1 in \S \ref{s-6}; 
\item[(iii)] means that locally the ($0$-)level curves 
of $H_i-H_j$, $i\neq j$, form a foliation 
by $1$-dimensional smooth curves of a small enough neighborhood 
of $p$. As (ii) above, this assumption will also be used in an essential way.
\end{itemize}

Let $K$ be the convex hull of the points $A_i(p)$, $i\in I(p)$, and denote by 
$\partial K$ its boundary, which is clearly an $i(p)$-gon. 
From Theorem \ref{main-r1} and its proof sketched in \S \ref{s-3} 
and completed in 
\S \ref{s-4} -- \S \ref{s-5} (see, in particular, Lemma \ref{extremepoints} in 
\S \ref{ss-prel} and Corollary \ref{Cor5} in 
\S \ref{ss-44}) we deduce the following:

\begin{Cor}\label{cor-ip-conv}
Assume all the hypotheses of Theorem \ref{main-r1} except conditions (i)--(ii)
and set $S(p)=\{i\in I(p)\mid A_i(p)\text{ is an extreme point of }K\}$. 
If $A_k(p)\notin \partial K$ for $k\in I(p)\setminus S(p)$ then 
the conclusion of Theorem \ref{main-r1} holds.
\end{Cor}

\begin{remark}\label{r-extreme}
In particular, Corollary \ref{cor-ip-conv} holds if $S(p)=I(p)$, i.e., all 
points $A_i(p)$, $i\in I(p)$, are extreme in $K$.
\end{remark}
 
We emphasize the fact that results similar to those above cannot hold 
for arbitrary $PA$ functions. Indeed, as we already noted, the 
requirement 
that $\partial \Phi/\partial {\bar z}\geq 0$ is crucial. In particular, 
it implies that the open sets
$\{M_i\}_{i=1}^{r}$ and the analytic functions $\{A_i(z)\}_{i=1}^{r}$ 
associated
with $\Phi$ have to be intimately related to each other. The latter statement 
is illustrated (and further reinforced) in the next example.

\begin{Ex}\label{ex-bla2}
Let $r=2$, $A_1(z)\equiv 1$ and $A_2(z)\equiv i$. Then the subharmonic function
$\varphi$ defined in Theorem~\ref{main-r1} becomes 
$\varphi(x,y)=\max(x,-y)$, that is,
$\varphi(x,y)=x$ if $x+y\ge 0$ and $\varphi(x,y)=-y$ for $x+y\le 0$. Hence
its derivative $\frac{2\partial \varphi}{\partial z}$ equals $1$ if $x+y\ge 0$ 
and $i$ for $x-y\le 0$, respectively. Theorem~\ref{main-r1} 
says (loosely) that among
all $PA$ functions $\Phi$ of the form $1\cdot \chi_{M_1}+i\cdot\chi_{M_2}$
for varying sets $M_1$ and $M_2$ (covering some neighborhood of the origin up
to a Lebesgue negligible set) $\frac{2\partial \varphi}{\partial z}$
is the only one that has a 
positive $\bar{z}$-derivative in the sense of distributions. To see why 
this is the case consider the following simple example: let $l$
be a line through the origin with unit normal $n=n_1+in_2$, so that  
$\bC\setminus l$ consists of two half-planes. Let $M_1$ be the one
with $n$ as interior normal to its boundary and $M_2$ the other half-plane.
Set $\Phi=1\cdot \chi_{M_1}+i\cdot\chi_{M_2}$. Then
$$\frac{\partial \Phi}{\partial \bar{z}}=\frac{1}{2}(1-i)(n_1+in_2)ds,$$
where $ds$ is Euclidean length measure along the common boundary $l$
to $M_1$ and $M_2$ (see Corollary~\ref{coro4}). Clearly, 
$\partial \Phi/\partial \bar{z}\ge 0$ only if 
$n_1+in_2=\frac{1}{\sqrt{2}}(1+i)$, i.e., if the line $l$ is given by 
$x+y=0$. In other words one must indeed have 
$\Phi=\frac{2\partial \varphi}{\partial z}$, where
$\varphi$ is the subharmonic function defined in Theorem~\ref{main-r1}. 
Note that in this particular 
example we used the fact that the 
boundaries of the $M_i$'s are $C^1$ in order to explicitly calculate the 
derivative of $\Phi$. Our theorems show that
the corresponding result is true in a much more general situation 
with no assumptions on the boundaries.
\end{Ex}
 
The local characterization of subharmonic functions with $PA$ 
derivatives is 
almost an immediate consequence of Theorem~\ref{main-r1} and shows that at 
generic points such functions are indeed maxima of a finite set of harmonic 
functions:

\begin{Cor}\label{loc-char} 
Let $\Psi$ be a subharmonic 
 function such that $\partial \Psi/\partial z$ is a $PA$ function with 
 decomposition given by \eqref{eq:MDA} and satisfying conditions (i)--(iii) of
Theorem~\ref{main-r1}. Then there exists a neighborhood $\widetilde N(p)$ of 
$p$ and harmonic functions $H_{i}$, $1\le i\le r$, defined in  
$\widetilde N(p)$ such that $\Psi(z)=\max_{1\le i\le r}H_{i}(z)$ a.e.~in 
$\widetilde N(p)$.
\end{Cor}

Let $\Phi\in\Si(U)$, so that by Notation \ref{not-Sigma} and 
\cite[Theorem 2.1.7]{H} the measure $\nu:=\partial \Phi/\partial {\bar z}$ is 
positive. 
%Assume that $p\in U$ satisfies conditions (i)--(iii) of 
%Theorem~\ref{main-r1} and 
Let further $p\in U$ and $N(p)$ be a neighborhood of $p$ such that 
$\overline{N(p)}\subset U$. Then the (positive) measure 
$\tilde\nu:=\chi_{_{\overline{N(p)}}}\cdot \nu$ extends to $\bC$ and 
there exists some analytic function 
$A$ such that $\Phi=C_{\tilde \nu}+A$ (as distributions) in $N(p)$, 
where $C_{\tilde \nu}$ is the 
Cauchy transform of $\tilde \nu$ defined by
$$C_{\tilde \nu}:=\frac{1}{\pi z}*\tilde \nu.$$
The above decomposition for $\Phi$ is a consequence of formula (4.4.2) in 
{\em op.~cit.~}asserting that $\Phi$ and $C_{\tilde \nu}$ have the same 
derivative 
with respect to
 $\partial/\partial {\bar z}$, so that by \cite[Theorem 4.4.1]{H} they must 
differ by an analytic function. Hence we also have the following corollary to 
Theorem~\ref{main-r1}.
 
\begin{Cor}\label{Cauchy} 
Let $\Phi\in\Si(U)$ be a $PA$ function with 
 decomposition given by~\eqref{eq:MDA} and set 
$\nu=\partial \Phi/\partial {\bar z}$. Assume that $p\in U$ 
satisfies conditions (i)--(iii) of Theorem~\ref{main-r1} and let 
$N(p)$ and $\tilde \nu$ be as above. Then 
$\Phi=C_{\tilde \nu}+A$ in $N(p)$, 
where $A$ is an analytic function and the positive measure 
$\tilde \nu$ is supported in a 
union of segments of level sets for the functions $H_i-H_j$, 
$1\le i\neq j\le r$. Moreover, $\nu$ may be locally described by means of its 
support in the sense of formula \eqref{eq:MDA3} 
(see Theorem~\ref{coro6} (3) in \S \ref{s-2}).
\end{Cor}

Note that the above results hold in a surprisingly great generality as
they assume no {\em \`a priori} knowledge of the differentiable structure of 
$\text{supp}\,\nu$. We construct an example 
showing that the picture is even more complex in non-generic cases and in 
particular that Corollary~\ref{loc-char} is not true if $p$ is special enough,
see Example~\ref{ex1} in \S \ref{s-6}.

The special case when the $A_{i}$ in Theorem~\ref{main-r1} are 
constant functions was treated in \cite{BR}. Our crucial 
Lemma \ref{lemma:decreasing} is mutatis mutandis generalized from that paper.
In the simpler situation 
of {\em loc.~cit.~}some additional global results were 
obtained. These show essentially that any (locally) $PH$ subharmonic 
function is globally (in $U$) a maximum of finitely many harmonic functions.
Example~\ref{ex1} in \S \ref{s-6} again shows that this is not true in general.
However, it is not difficult to get complete results in the case 
when only two functions are involved, see \S \ref{s-2}. It would be 
interesting to establish when a subharmonic function with a $PA$ 
derivative is globally a maximum of finitely many 
harmonic functions (cf.~Problem~\ref{pb2} in \S \ref{s-6}). 

\section{Derivatives of Sums}\label{s-2}
%: Regular Boundary Cases

Recall the canonical piecewise decomposition of a $PH$ function from 
\S \ref{s-12} 
(cf.~Definition \ref{d-r2} with $X=H(U)$). If $\Psi(z)$ is a 
$PH$ subharmonic function of the form 
\eqref{eq:MDH} then the support of the associated Riesz measure $\De \Psi$ 
equals $\Gamma:=U\setminus \bigcup_{i=1}^{r}M_{i}$. Indeed, it is clear that 
$\text{supp}(\De \Psi)\subseteq \Gamma$. For the reverse inclusion note that 
$\Psi$ is harmonic in a neighborhood of any point 
$p\in\Gamma\setminus \text{supp}(\De \Psi)$. If such a point exists 
one can find $i\neq j$ so that any neighborhood 
of $p$ intersects $M_i$ and $M_j$, and then $H_i$ and $H_j$ both agree with 
$\Psi$ in some neighborhood of $p$ hence $H_i=H_j$ (by the unique continuation 
property), which is a contradiction.

In this section we first discuss the case of a $PA$ function $\Phi$ 
with canonical 
piecewise decomposition as in Definition \ref{d-r2} such that the 
corresponding set $\Ga=U\setminus \bigcup_{i=1}^{r}M_{i}$ is a locally finite 
union of piecewise $C^1$-curves. We show that if the distribution derivative 
$\pa \Phi/\pa \bar{z}$ is positive then this measure is determined in a 
simple way by its support, see Theorem~\ref{coro6} (3) below. 
Note that in view of Lemma \ref{l-can} (v) 
a situation where $\Ga$ is piecewise smooth 
occurs if one considers a $PA$ function of the form 
$\Phi=\sum_{1\le i\le r}(\pa H_i/\pa z)\chi_i$, where 
$\Psi=\sum_{1\le i\le r}H_i\chi_i$ is a continuous $PH$ function (for 
instance, $\Psi$ could be the maximum of finitely many harmonic functions). 
In this case we show that the continuity assumption implies that 
$\Phi$ is actually the distribution derivative of $\Psi$ (without any
$C^1$-assumptions on $\Ga$).

% which occurs for instance if $\Psi(z)$ is  We show that if 
%this assumption is satisfied then $\Gamma$ is essentially a 
%union of level 
%curves $\{z\in U\mid H_{i}(z)=H_{j}(z)\}$, $1\le i\neq j\le r$, and that the 
%measure $\De \Psi$ is determined by its support in the sense of formula 
%\eqref{eq:MDA3} below (Proposition~\ref{coro5} and ). 
%As it will become transparent from the proofs, these results hold actually 
%even under weaker regularity conditions such as Lipschitz continuity for the 
%curves in $\Ga$.

We start with the case when only two functions are involved. Assume that 
$\Phi(z)$ is defined in a domain $U$ and that
there exists a smooth curve $\Gamma\subset U$ dividing $U$ into two open
connected components $U=M_{1}\cup\Gamma\cup M_{2}$ such that
$\Phi(z)=A_{i}(z)$ in $M_{i}$, $i=1,2$, where $A_{i}(z)$ is a
function analytic in some neighborhood of $M_{i}$. In particular, $\Phi(z)$ is 
a $PA$ function. 

\begin{Lemma}\label{lm:tang} 
If $\nu:=\partial \Phi(z)/\partial \bar z\ge 0$ in the sense
of distribution theory (i.e., $\nu$ is a positive measure) then at each
point $\tilde z$ of $\Gamma$ the
tangent line $l(\tilde z)$ to $\Gamma$ is orthogonal to
$\overline{A_{1}(\tilde z)}-\overline{A_{2}(\tilde z)}$ and the measure 
$\nu$ at $\tilde z$ equals
$$\frac{\vert
A_{1}(\tilde z)-A_{2}(\tilde z)\vert ds}{2},$$
where $ds$ denotes length measure along $\Gamma$.
\end{Lemma}

Lemma~\ref{lm:tang} is an immediate consequence of the following well-known 
result, see e.g.~\cite[Theorem 3.1.9]{H}.

\begin{Prop}\label{propo3}
Let $Y\subset  X$ be open subsets of $\bR^k$ such that $Y$ has
             a $C^1$-boundary $\partial Y$ in $X$ and let $u\in C^1(X)$.
If $\chi_{_Y}$ denotes the
characteristic function of $Y$, $dS$ the Euclidean surface measure
on $\partial Y$
and $n$ the interior unit normal to $\partial Y$ then
$$\partial_j (u\chi_{_Y})=(\partial_j u) \chi_{_Y}+un_jdS,$$
where $\partial_j$ and $n_{j}$ are the partial derivative with respect to 
the $j$-th coordinate and the $j$-th component of $n$, respectively.
\end{Prop}

\begin{Cor}\label{coro4}
In the notation of Proposition~\ref{propo3} one has
\begin{equation}\label{eq:barz}
\begin{split}
&\frac{\partial (u\chi_{_Y})}{\partial\bar z}=\left(\frac{\partial
u}{\partial \bar z}\right)\!\chi_{_Y}+\frac{1}{2}u(n_1+in_2)ds,\\
&\frac{\partial (u\chi_{_Y})}{\partial z}=\left(\frac{\partial
u}{\partial z}\right)\!\chi_{_Y}+\frac{1}{2}u(n_1-in_2)ds.
\end{split}
\end{equation}
\end{Cor}

\begin{proof}[Proof of Lemma \ref{lm:tang}]
Suppose that the function $\Phi(z)=A_{1}(z)\chi_{1}(z)+A_{2}(z)\chi_{2}(z)$ 
satisfies the conditions of Lemma \ref{lm:tang}, where $\chi_{i}$ is the 
characteristic function of $M_{i}$, $i=1,2$. Corollary~\ref{coro4} implies in 
particular that $\nu$ is supported on the smooth separation curve $\Gamma$
and that with an appropriate choice of co-orientation one has 
$\nu=\frac {(A_{1}-A_{2})nds}{2}$, which proves the lemma.
\end{proof}

Proposition~\ref{propo3} remains true if the boundary of
$Y$ is assumed to be only 
piecewise $C^1$ or just Lipschitz continuous (cf.~{\em op.~cit.}). We may 
therefore apply it to functions of the form  
$$\max_{1\le i\le r} H_{i}(z)=\sum_{i=1}^{r}H_{i}(z)\chi_{i}(z)$$
in $U$ and get the 
 description of their derivatives given in the introduction. In this 
 case the normal $n$ is defined a.e.~with respect to length measure on the 
 boundary and the equality in Corollary~\ref{loc-char} is interpreted in this 
sense. 

\begin{notation}\label{not-1}
Given a $PH$ function $\Psi(z)=\sum_{i=1}^{r}H_{i}(z)\chi_{i}(z)$ as in 
\eqref{eq:MDH} let $\Gamma_{\Psi}= 
 U\setminus \bigcup_{i=1}^{r}M_{i}$ and denote by $\Gamma_{\Psi}^d$ the set 
of points where the normal to $\Gamma_{\Psi}$ is not defined. In similar 
fashion, for a $PA$ function
$\Phi(z)=\sum_{i=1}^{r}A_{i}(z)\chi_{i}(z)$ as in \eqref{eq:MDA} 
we set $\Gamma_{\Phi}= 
 U\setminus \bigcup_{i=1}^{r}M_{i}$ and let $\Gamma_{\Phi}^d$ be the set 
of points where the normal to $\Gamma_{\Phi}$ is not defined.
\end{notation}

Essentially 
the same arguments yield the following generalization of Lemma \ref{lm:tang}.

\begin{Theorem}\label{coro6}
Let 
$$\Phi(z)=\sum_{i=1}^{r}A_{i}(z)\chi_{i}(z)$$
be a $PA$ function in a simply connected domain $U\subset \bC$ such that
\begin{itemize}
\item[(i)] $\Gamma_{\Phi}$ is a locally finite union of piecewise 
$C^1$-curves;
%\item[(ii)] 
%$\Gamma_{\Phi}^{d}$ has measure $0$ with respect to length measure 
%$ds$ on $\Gamma_{\Phi}$;
\item[(ii)] $\partial \Phi/\partial\bar z\geq 0$.
\end{itemize}
Let $H_{i}$, $1\le i\le r$, be real-valued harmonic functions as in 
\eqref{a-h}.
%a primitive function of $A_{i}$.
Then for any $z\in \Gamma_{\Phi}\setminus \Gamma_{\Phi}^{d}$ there is a 
neighborhood $N(z)$ such that 
\begin{enumerate}
\item[(1)] $N(z)\setminus \Ga_\Phi$ consists of two components 
$N(z)_i$, $N(z)_j$ such that $\Phi(z)=A_k(z)$ in $N(z)_k$ for $k=i,j$; 
\item[(2)] $N(z)\cap \Ga_\Phi$ is contained in a level curve of $H_i-H_j$
%$\Ga_{ij}:=\{z\in U\mid H_{i}(z)=H_{j}(z)\}$; 
for some $i,j$;
\item[(3)] In $N(z)$ one has
\begin{equation}
\partial \Phi(z)/\partial\bar  z=\frac{\vert A_{i}(z)-A_{j}(z)\vert ds}{2}. 
\label{eq:MDA3}
\end{equation}
The restriction of $\partial \Phi(z)/\partial\bar  z$ to $U\setminus \Gamma_{\Phi}^{d}$, determined locally by \eqref{eq:MDA3}, 
%The right-hand side of \eqref{eq:MDA3} 
extends to a measure $\mu$ on $U$
%on $\Ga_\Phi$ 
which is absolutely continuous with respect to length measure on $\Ga_\Phi$. Furthermore 
%and such that \eqref{eq:MDA3} (for suitabholds in 
$\partial \Phi(z)/\partial\bar  z=\mu$ in
$U$. Moreover, if any two level 
curves $\Ga_{ij}$, $\Ga_{kl}$ with $i<j$, $k<l$, 
$(i,j)\neq (k,l)$ intersect in at most 
a finite number of points, then the measure $\mu$ hence also 
$\partial \Phi(z)/\partial\bar  z$ is determined by its support $\Ga_\Phi$.
\end{enumerate}
\end{Theorem}

%$\Gamma_{\Phi}\setminus \Gamma_{\Phi}^{d}$ is a union of pieces of the 
%level curves 
%$\Ga_{ij}:=\{z\in U\mid H_{i}(z)=H_{j}(z)\}$, $1\le i\neq j\le r$, 
%and on such a level curve one has
%Then
%any smooth part of 
%$\Ga_\Phi$ belongs to a unique $\Ga_{ij}$. Moreover, the measure $\mu$ on 
%the whole $\Ga_\Phi$ given by \eqref{eq:MDA3}. Then 

\begin{proof}
Assertions (1), (2) and identity \eqref{eq:MDA3} are direct 
consequences of Lemma~\ref{lm:tang}.
Since by (i) $\Gamma_{\Phi}$ is a locally finite union 
of piecewise $C^1$-curves 
the set $\Gamma_{\Phi}^{d}$ has measure $0$ with respect to length measure 
$ds$ on $\Gamma_{\Phi}$ and thus the measure $\mu$ extending the 
right-hand side of \eqref{eq:MDA3} to $\Gamma_{\Phi}$ exists. It remains to 
show that 
\begin{equation}\label{ex-cauchy}
\partial \Phi/\partial\bar  z=\mu.
\end{equation} 
Note that 
$\partial \Phi/\partial\bar  z=\mu+G$, where $G$ is a sum of Dirac 
measures supported at (singular) points in $\Gamma_{\Phi}^{d}$. Consider now 
a singular 
point $p\in \Gamma_{\Phi}^{d}$, a small neighborhood $N$ of $p$, and the 
Cauchy transform $C_{\tilde \mu}$ 
of (the extension to $\bC$ of) the measure 
$\tilde \mu:=\chi_{_{\overline{N}}}\cdot \mu$. Suppose that locally at $p$ 
the measure $G$ is given by $c\delta_p$ for some $c\ge 0$. Then the function 
$$
\Phi-C_{\tilde\mu}-\frac{c}{z-p}
$$
is analytic at $p$. On the other hand, $\Phi$ is bounded and by the 
classical Plemelj-Sokhotski formulas (cf., e.g., \cite[\S 3.6]{BG}) the Cauchy 
transform $C_{\tilde\mu}$ has at most a logarithmic singularity at $p$. It 
follows that $c=0$, which proves \eqref{ex-cauchy}. For the last statement in 
part (3) of the theorem note that the assumption on the level curves made 
there guarantees that each regular point of $\Ga_\Phi$ belongs to a unique 
$\Ga_{ij}$, hence in view of \eqref{eq:MDA3} the measure 
$\partial \Phi/\partial\bar z$ is locally determined by $\Ga_{ij}$.
%Using the second relation in~\eqref{eq:barz} it 
%suffices to 
%check locally at a point $z\in \Gamma\setminus \Ga_{\Psi}^d$ where 
%$H_{i}(z)=H_{j}(z)$ for some $i\neq j$ that the contributions 
%$\frac{1}{2}H_{i}(n_1-in_2)ds$ and $\frac{1}{2}H_{j}(n_1-in_2)ds$ 
%cancel; this follows since 
%the normal occurs twice with opposite directions. Note that if $\Psi$ is 
%subharmonic then Lemma \ref{lm:tang} implies that $\Gamma_{\Psi}$ consists of 
%pieces of level curves to the functions $H_{i}-H_{j}$, $1\le i\neq j\le r$. 
%As $\Psi$ is continuous these level curves have to be given by 
%$\{z\in U\mid H_{i}(z)=H_{j}(z)\}$.
\end{proof}

In the remainder of this paper we will see that results similar to 
Theorem \ref{coro6} actually  
hold without local regularity assumptions as in 
%of Proposition~\ref{coro5} 
Theorem~\ref{coro6} (i).

Obviously, a $PH$ function $\Psi$ has a $PA$ derivative almost everywhere.
However, this is not necessarily the same as the 
distribution derivative of $\Psi$. The next result shows that this is true 
for continuous $PH$ functions.

\begin{Prop}\label{coro5}
If the canonically decomposed PH function 
$$\Psi(z)=\sum_{i=1}^{r}H_{i}(z)\chi_{i}(z)$$ is continuous in $U$ 
(cf.~\S \ref{s-12})  
%and $\Gamma_{\Psi}$ is a locally finite union of piecewise $C^1$-curves 
then 
\begin{equation}
    \partial \Psi(z)/\partial z=
\sum_{i=1}^{r}A_{i}(z)\chi_{i}(z)
 \label{eq:MDA2}
\end{equation}
in the sense of distributions, where 
$A_{i}:=\partial H_{i}/\partial z$, $1\le i\le r$. 
%Moreover, if $\Psi(z)$ is 
%subharmonic then  
%$\Gamma_{\Psi}$ consists of pieces of the level curves 
%$\{z\in U\mid H_{i}(z)=H_{j}(z)\}$, $1\le i\neq j\le r$.  
\end{Prop}  

\begin{proof}
Let $\Ga_\Psi$ be as in Notation \ref{not-1}. By Lemma \ref{l-can} (5) 
$\Ga_\Psi$ is contained in the zero set of the function 
$g=\prod_{1\le i<j\le r}(H_i-H_j)$. Let $p\in\Ga_\Psi\setminus \Ga_\Psi^d$ be 
a regular point of $\Ga_\Psi$ and $N$ be a small (open) neighborhood of $p$. 
Let further $N_{\pm}$ be $N$ intersected with the two sides of $\Ga_\Psi$. 
It follows that 
$N_+\subset M_i$ and $N_{-}\subset M_j$ for some $i\neq j$ if $N$ is 
small enough, and the restriction of $\Psi$ to $N$ is a smooth function plus 
$f\chi_i$, where $f\equiv 0$ in $\Ga_\Psi$. Then $\pa (f\chi_i)/\pa z$ is a 
function in $N$ and we conclude that 
$\pa \Psi/\pa z=\sum_{i=1}^{r}A_i\chi_i+G$, where $G$ is a distribution 
supported at the singular points $\Ga_\Psi^d\subset \Ga_\Psi$. Since 
$\Ga_\Psi^d$ is a discrete set, by choosing a continuous solution $h$ to 
$\pa h/\pa z=\sum_{i=1}^{r}A_i\chi_i$ we get a continuous solution $\Psi-h$ 
to $\pa(\Psi-h)/\pa z=G$ and it follows that $G\equiv 0$, which proves the 
proposition.
\end{proof}

%The above results show in particular that it is rather easy to understand
%the Riesz measure of a continuous 
%subharmonic $PH$ function if one knows {\em \`a priori} 
%that its support is a locally finite union of 
%piecewise $C^1$-curves.

% which is not really 
%surprising. Indeed, consider the example of just one function $\Psi=H$ 
%a.e.~in $U$. 
%Then standard theorems (see~\cite{H}) imply that $\Psi$ is in fact 
%harmonic everywhere in $U$. 
%It is reasonable to expect that there should 
%be a direct general method to check that the support of $\nu$ is indeed at 
%least a locally finite countable union of 
%$C^1$-curves. An apparatus for proving results in this direction 
%seems to be contained in \cite{Mat}.

\section{Local Characterization in Generic Cases: Sketch of Proof}\label{s-3}

In this section we give an equivalent formulation of Theorem~\ref{main-r1} 
and sketch its proof. Under some mild 
non-degeneracy assumptions, this 
provides a local description of functions with positive (distributional) 
$\bar z$-derivative which is equal a.e.~to one of a finite 
number of given analytic functions. 

Let us first fix notations and assumptions.

\begin{notation}\label{not-2}
Let $\{M_{i}\}_{i=1}^r$, $r\ge 2$, be a finite family of disjoint 
open subsets of a simply connected domain $U\subset{\bC}$ covering $U$
up to a set of zero Lebesgue measure and denote by
$\chi_{i}$ the characteristic function of $M_{i}$. Given a family 
$\{A_i(z)\}_{i=1}^r$ of pairwise distinct analytic functions in $U$ define 
the (measurable) 
function 
$$\Psi(z)=\sum_{i=1}^{r}A_{i}(z)\chi_{i}(z).$$ 
Fix a point $p\in U$. As in \eqref{a-h} we let 
$$H_{i}(z)=\Re\left[\int_{p}^z A_{i}(w)dw\right], \quad 1\le i\le r.$$ 
Note that each $H_i$ is a well-defined harmonic function in 
$U$ satisfying $ \partial H_{i} /\partial z=\frac{1}{2}A_{i}(z)$.
If $r\ge 3$ we associate to each triple $(i,j,k)$ of distinct indices in 
$\{1,\ldots,r\}$ the following ``critical set''
\begin{equation}\label{eq-crit-set}
\Gamma_{i,j,k}=\{z\in U\mid\text{$A_i(z)$, $A_j(z)$, $A_k(z)$ are collinear}\}.
\end{equation}
Alternatively, $\Gamma_{i,j,k}$ consists of the set of $z\in U$ such that 
$A_{i}(z)-A_{k}(z)$ and $A_{j}(z)-A_{k}(z)$ are linearily dependent over the reals. 
This is the set where the gradients of $H_{i}-H_{k}$ 
and $H_{j}-H_{k}$ are parallel, or 
equivalently, the level curves through $z$ to these functions are 
parallel. Clearly, 
$\Gamma_{i,j,k}$ is either a real analytic curve or else 
there exists $c\in\bR$ such that $A_{i}(z)-A_{k}(z)=c(A_{j}(z)-A_{k}(z))$
for all $z\in U$.
\end{notation}

%The following lemma is left as an (elementary) exercise.

%\begin{Lemma}\label{s3}
%Let $\fS_3$ be the symmetric group on $3$ elements. For each triple $(i,j,k)$ 
%of distinct indices in 
%$\{1,\ldots,r\}$ and $\pi\in\fS_3$ one has 
%$\Ga_{i,j,k}=\Ga_{\pi(i),\pi(j),\pi(k)}$.
%\end{Lemma}

In this notation Theorem \ref{main-r1} may then be restated as follows. 
Suppose -- using the labeling in the theorem -- that $i(p)=r$ 
(cf.~\eqref{ip}), assume that
$\partial \Psi/\partial {\bar z}\geq 0$ as a distribution
supported in $U$ and let $p\in U$ be such that 
%any neighborhood $N$ of $p$ intersects each $M_{i}$, 
%in a set of a positive Lebesgue measure. 
\begin{itemize}
\item[(i)] $p\in\overline{M}_i$, $1\le i\le r$;
\item[(ii)] There is no critical set $\Gamma_{i,j,k}$ that contains $p$;
\item[(iii)] $A_{i}(p)\neq A_{j}(p)$ for $1\le i\neq j\le r$, i.e., 
$p$ is a non-singular point of $H_{i}-H_{j}$.
\end{itemize}
Then there exists a neighborhood $\widetilde N(p)$ of $p$ such that
$$\Psi=2\partial \varphi /\partial z\text{ a.e.~in } \widetilde N(p),$$
where $\varphi$ is the subharmonic function defined by
$$ \varphi(z)=\max_{1\le i\le r}H_i(z).$$

\begin{remark}\label{rem-1}
Generically, the sets $\Gamma_{i,j,k}$ are curves and 
so conditions (ii) and (iii) above hold outside some real analytic set.
\end{remark}

\noindent
{\bf Strategy of the proof and two fundamental lemmas.} 
The proof of Theorem~\ref{main-r1} is rather technical and the main parts 
of the argument are contained in Lemma~\ref{lemma:decreasing} 
and Lemma~\ref{spaghetti} below, which to some extent hold independently 
of condition (ii) in  Theorem~\ref{main-r1}. We will now show that 
Theorem~\ref{main-r1} follows in fact from these two lemmas. 
First, a convenient reformulation of the 
conclusion of Theorem~\ref{main-r1} is that for $1\le i\le r$ one has 
$\chi_{i}=1$ a.e.~in the set where $\varphi(z)=H_{i}(z)$, and this is what we 
will actually show. Clearly, it is enough to prove this statement for 
$i=1$. 

\medskip

\noindent
{{\bf Assumption I.}} 
{\em By considering the function $\Psi - A_{1}$ and using the fact that 
$A_{1}$ is analytic in $U$ (hence 
$\partial A_{1}/\partial {\bar z}=0$) we may assume without loss of 
generality that} 
\begin{equation*}
A_{1}(z)=H_{1}(z)=0\text{ {\em for} }z\in U,\tag*{(I)}
\end{equation*}
{\em which we do, except when otherwise stated, throughout the remainder 
of this section as well as in \S \ref{s-4} and \S \ref{s-5}}.

\medskip 

Define now 
\begin{equation}\label{w-set}
\begin{split}
&W=W_{1}(p):=\left\{z\in U\mid \varphi(z)=0
\right\},\\
&W_i(p):=\left\{z\in U\mid \vf(z)=H_i(z)
\right\},\,2\le i\le r.
\end{split}
\end{equation}
We have to prove that $\Psi=0$ a.e.~in
$N\cap W$, or equivalently $\Psi=0$ a.e.~in
$N\cap \mathring{W}$ for some small enough neighborhood $N$ of $p$, 
where $\mathring{W}$ denotes the interior of $W$. 

The first lemma asserts that 
    $\chi_{1}$ is 
    increasing along every path along which all functions 
$H_{i}$, $2\le i\le r$,  
    are decreasing.

\begin{Lemma}\label{lemma:decreasing}
Let $p\in U$ satisfy all the assumptions of Theorem~\ref{main-r1} except 
condition (ii). If $\gamma$ is a piecewise $C^1$-path from $z_1=\ga(0)$ to 
$z_2=\ga(1)$ such that
that each of the functions 
$[0,1]\ni t\mapsto H_{i}(\gamma(t))$, $2\le i\le r$, is decreasing then 
\begin{equation}\label{chi:test}
(\chi_{1}*\phi)(z_{1})\leq (\chi_{1}*\phi)(z_{2})
\end{equation}
for any positive test function $\phi$ with $\text{{\em supp}}\,\phi$ small 
enough.
\end{Lemma}   

The second lemma guarantees that enough many points may be 
reached by paths of the form given in Lemma~\ref{lemma:decreasing}.
To make a precise statement we need the following definition: to each 
$z\in U$ we associate the set
$$V(z)=\{\ze\in U\mid \exists\text{ piecewise $C^1$-path from $z$ to $\ze$ 
along which all $H_{i}$ decrease}\}.$$

\begin{definition}\label{def-lim}
Given $p\in U$ and two subsets $M,X\subset U$ with $p\in\overline{M}$ 
we say that $V(z)$ {\em tends to}
$X$ {\em through} $M$ {\em as} $z\to p$, which we denote by 
$\lim_{M\ni z\to p}V(z)=X,$
if for each $\al\in X$ and any sequence $\{z_n\}_{n\in \bN}\subset M$ 
converging to $p$ 
one has $\al\in V(z_n)$ for all but finitely many indices $n\in \bN$.
\end{definition}

\begin{Lemma}\label{spaghetti}
Let $p\in U$ satisfy all the assumptions of Theorem~\ref{main-r1}, in 
particular $p\notin \Gamma_{i,j,1}$ for any $i,j$. Then there is a 
neighborhood $N$ of $p$ with 
$$
\lim_{U\ni z\to p}V(z)=N \cap \mathring{W}.
$$
\end{Lemma}

\begin{remark}
Note that there are actually no sets $\Ga_{i,j,k}$ at all if $r=2$ in 
Lemma~\ref{spaghetti}.  
\end{remark}

\begin{proof}[Theorem~\ref{main-r1}: outline of the proof]  
As noted in the paragraph preceding Lemma~\ref{lemma:decreasing}, we have to 
show that there exists a sufficiently small neighborhood $N$ of $p$ such that 
$\Psi=0$ a.e.~in $N\cap \mathring{W}$. This is trivially true 
if $W$ has no interior points (i.e., if $\mathring{W}$ has zero Lebesgue 
measure) 
and so we may assume that $\mathring{W}$ has positive Lebesgue measure. 

Let now $\{\phi_{s}\}_{s\in\bN}$ be a 
sequence of test 
functions satisfying $\text{supp}\,\phi_{s}\to \{0\}$ as $s\to \infty$ 
and $\int\!\phi_{s}d\lambda =1$, $s\in\bN$, where $\lambda$ denotes 
Lebesgue measure. Note that $\{\phi_{s}* \chi_{1}\}_{s\in\bN}$ converges in
$L^1_{loc}$ to $\chi_{1}$. In particular, this implies that for all 
$\epsilon>0$, $\delta>0$ there exist a sufficiently large $s(\eps,\de)\in\bN$ 
such that if $s\in\bN$, $s\ge s(\eps,\de)$, there is 
a point $z_{1}=z_1(\eps,\de,s)\in U$ satisfying 
\begin{equation}\label{eq:tech}
\vert z_{1}-p\vert<\delta
\text{ and }(\phi_{s}* \chi_{1})(z_{1})>1-\epsilon.
\end{equation}
To see this let $N_\delta=\{z\in U\mid\vert z-p\vert<\delta\}$ and suppose 
that $(\phi_{s_k}* \chi_{1})(z)\le 1-\eps$ for some infinite sequence 
$\{s_k\}_{k\in\bN}$ and almost all $z\in N_{\de}$. Then
$$\int_{N_\delta}\left\vert (\phi_{s_k}*\chi_{1})(z)-\chi_{1}(z)\right\vert 
d\lambda(z)>\epsilon\lambda(M_1\cap N_\delta)$$
and since by assumption $\lambda(M_1\cap N_\delta)>0$ this contradicts the 
fact that $\{\phi_{s_k}* \chi_{1}\}_{s\in\bN}$ converges to $\chi_{1}$ in
$L^1_{loc}$ as $k\to\infty$, so that~\eqref{eq:tech} must hold.

From \eqref{chi:test} and~\eqref{eq:tech} it follows that
$(\phi_{s}* \chi_{1})(z)>1-\epsilon$ for $z\in V(z_{1})$, which together 
with the identity $\phi_{s}* 1=1$ yields $\big(\phi_{s}*
\sum_{i=2}^r\chi_{i}\big)(z)< \epsilon$ and therefore
\begin{equation}\label{eq:convol}
\begin{split}
\left|(\phi_{s}* \Psi)(z)\right|&=\left|\int\!\phi_{s}(z-\zeta)
\Psi(\ze)d\lambda(\zeta)\right|\\
&\leq \epsilon \max_{2\le d\le r}\,
\sup_{\ze\in z-\text{supp}\,\phi_s}|A_{d}(\ze)|=:\epsilon C_s(z),\quad 
z\in V(z_{1}).
\end{split}
\end{equation}

Now we assume in addition that all the conditions of Theorem~\ref{main-r1}
and Lemma~\ref{spaghetti} are true. 
Fix $\eps>0$. The arguments 
above show that one can construct a sequence 
$\{z_n\}_{n\in\bN}\subset U$ such that
\begin{equation}
|z_n-p|<\frac{1}{n}\text{ and }(\phi_{s_n}* \chi_1)(z_n)>1-\eps
\label{eq:foeljd}
\end{equation}
for some strictly increasing sequence of positive integers $\{s_n\}_{n\in\bN}$.
By Lemma~\ref{spaghetti} there exists a neighborhood $N$ of $p$ such that 
each $z\in N \cap \mathring{W}$ belongs to all but finitely many sets 
$V(z_{n})$, $n\in \bN$. 
Combined with~\eqref{eq:convol} this shows that for every 
$z\in N \cap \mathring{W}$ 
there exists $n_z\in\bN$ such that
\begin{equation}\label{eq:ineq}
\left|(\phi_{s_n}* \Psi)(z)\right|\le C_{s_n}(z)
\eps\text{ for }n\ge n_z.
\end{equation}
Since $A_d$, $2\le d\le r$, are analytic functions and 
$\text{supp}\,\phi_{s_n}\to \{0\}$, $n\to \infty$, it follows 
from~\eqref{eq:convol} that by shrinking the neighborhood $N$ (if necessary) 
one can find $C>0$ such that $C_{s_n}(z)\le C$ for  
$n\in\bN$ and $z\in N \cap \mathring{W}$. Together with~\eqref{eq:ineq} and 
the fact 
that $\lim_{n\to\infty}\phi_{s_n}* \Psi=\Psi$ in $L^1_{loc}$
this clearly implies that $\Psi=0$ a.e.~in $N\cap \mathring{W}$, which 
proves Theorem~\ref{main-r1}.
\end{proof}

\section{Proof of Lemma~\ref{spaghetti}}\label{s-4}

To complete the proof of Theorem~\ref{main-r1} it remains to 
show Lemma~\ref{lemma:decreasing} and Lemma~\ref{spaghetti}. We start with 
the latter, which we prove in this section.

\subsection{Preliminaries}\label{ss-prel}

Let $A(z)$ be an analytic function defined in a neighborhood of some point 
$z_0\in\bC$ and set $H(z):=\Re\!\left[\int_{z_{0}}^z A(w)dw\right]$, so 
that $\partial H(z) /\partial z=\frac{1}{2}A(z).$
The directional derivative of $H$ with respect to a 
 complex number $v=\alpha+\beta i$ is given by
 \begin{equation}
  D_{v} H(z)=\alpha\partial H(z)/\partial x+\beta \partial 
 H(z)/\partial y= \Re\left[vA(z)\right]
     \label{eq:riktnder}
 \end{equation}
and the gradient of $H(x,y)$ 
considered as a vector in $\bC$ is just 
\begin{equation}\label{eq-grad}
\na H(x,y)=2\partial 
H(z)/\partial \bar z=\overline{A(z)}.
\end{equation}
If $A(z_{0})\neq 0$ 
then $z_{0}$ is a non-critical point for $H(z)$ and locally the $0$-level 
curves of $H$ form a foliation 
by $1$-dimensional smooth curves of a small enough neighborhood $N$ 
of $z_{0}$ (\cite[Theorem 5.7]{Spivak}).
In particular, the ($0$-)level curve $C_{H}$ of $H$ through $z_{0}$ 
divides $N$ into two
components
$$N_{H}^+=\{z\in N\mid H(z)>0\},\quad N_{H}=\{z\in N\mid H(z)<0\}.$$
Correspondingly, 
the tangent to $C_{H}$ at $z_0$ divides 
the plane into two opposite half-planes 
\begin{equation*}
\begin{split}
&\tau(z_0)^+=\left\{v+z_0\mid v\cdot \na H(z_0)\geq 0\}=
\{v+z_0\mid \Re\left[vA(z_0)\right]\geq 0\right\},\\
&\tau(z_0)=\left\{v+z_0\mid v\cdot \na H(z_0)\le 0\}=
\{v+z_0\mid \Re\left[vA(z_0)\right]\le 0\right\}.
\end{split}
\end{equation*}

We now return to the functions $A_i$, $1\le i\le r$, suspending for the 
moment Assumption I in \S \ref{s-3} stating that $A_1=0$. As before, we 
suppose that  $A_i(p)\neq A_j(p)$ if $i\neq j$. Consider the convex hull 
$K$ of the points $A_i(p)$, $1\le i\le r$. For each $i$ define the dual cone 
(with vertex at $p$) to the sector consisting of all rays from 
$\na H_i(p)=\overline{A_i(p)}$ to points in the complex dual $\overline K$  
%bounded by the lines 
%supported by two edges of $K$ meeting at 
%$A_i(p)$ -- actually, this is dual to the corresponding sector with 
%vertex $\overline A_i(p)$ in the complex dual $\overline K$ -- 
by 
\begin{equation}\label{eq-cone}
\begin{split}
\sigma_i(p):&=\bigcap_{k\in K}
 \left\{v+p\mid v\cdot (\bar k-\na H_i(p))\le 0\right\}\\
&=\bigcap_{j\neq i}^{r}
 \left\{v+p\mid v\cdot (\na H_j(p)-\na H_i(p))\le 0\right\}\\
&=\bigcap_{j\neq i}^{r} 
\left\{v+p\mid \Re\left[v(A_j(p)-A_i(p))\right]\le 0\right\}.
\end{split}
\end{equation}
Clearly, this cone is the infinitesimal analogue of the set $W_i(p)$ 
defined in (\ref{w-set}). The interior of $\sigma_i(p)$ contains the 
directions in which $H_i$ grows faster (up to the first order) than any other 
$H_k$, $k\neq i$.

There are several possibilities for the cone $\sigma_i(p)$: (a) it may have a 
top angle strictly between $0$ and $\pi$, in which case we say that it is a 
{\em pointed cone} (b) it consists just of the point $p$ or (c) it is either 
a line, a half-line or a half-plane.

The next lemma is a direct consequence of basic convex geometry.

\begin{Lemma}\label{convexhull}
With the above notations and assumptions the following holds: 
\begin{itemize}
\item[(i)] If $A_i(p)$ lies in the interior of $K$ then $\sigma_i(p)=\{p\}$;
%\item[(ii)]  $p\in \Gamma_{i,j,k}$ if and only if $A_i,A_j,A_k$ are collinear;
\item[(ii)] If $K$ is not a segment then $A_i(p)$ is an extreme point of $K$ 
if and only if $\sigma_i(p)$ is a pointed cone.
\end{itemize}
\end{Lemma}

Now consider condition (ii) in Theorem \ref{main-r1}, which is also part of 
the assumptions of Lemma \ref{spaghetti}. By Lemma \ref{convexhull} (ii) 
this condition is strictly stronger than the hypothesis in the following lemma.

\begin{Lemma}\label{extremepoints} 
Assume that the only points $A_i(p)$ contained in the boundary $\partial K$ 
of $K$ are extreme points. If $S(p)=\{i\in\{1,\ldots,r\}\mid A_i(p) \text{ is 
an extreme point of }K\}$ then:
\begin{itemize}
\item[(i)] $\max_{1\le i\le r}H_i(z)=\max_{i\in S(p)}H_i(z)$ in a neighborhood 
of $p$; 
\item[(ii)] There is a neighborhood $N$ of $p$ such that
$\bigcup_{i\in S(p)} N\cap W_i=N$.
\end{itemize}
\end{Lemma}

\begin{proof}
Clearly, (ii) follows from (i). Let now $j\notin S(p)$, so that by 
Lemma \ref{convexhull} and the assumption of Lemma \ref{extremepoints} one has 
$\sigma_j(p)=\{p\}$. 
This means that for each ray from $p$ in the unit vector direction $v\in S^1$ 
there is at least one $H_i$, $i\in S(p)$, such that 
$$
v\notin \left\{u+p\mid u\cdot (\na H_i(p)-\na H_j(p))\le 0\right\}.
$$ 
Thus, for each $v\in S^1$ there is a product neighborhood 
$I(v)\times J(v,p)\subset S^1\times U$ of $\{v\}\times \{p\}$ such that 
there exists $i=i(v)\in S(p)$ so that the continuous function 
$u\cdot (\na H_i(z)-\na H_j(z))$ is positive if $(u,z)\in I(v)\times J(v,p)$. 
By the compactness of $S^1\times \{p\}\subset S^1\times U$, a finite number 
of neighborhoods 
$I(v_l)\times J(v_l,p)$, $1\le l\le s$, cover $S^1\times \{p\}$. 
Hence the neighborhood $J(p):=\bigcap_{1\le l\le s} J(v_l,p)$ of $p$ has the 
property that
along each ray from $p$ with direction $v\in S^1$ there is some 
$i\in S(p)$ such that 
$H_i(z)>H_j(z)$ if $z\in J(p)\setminus \{p\}$, which proves (i).
\end{proof}

For the rest of this section we may again (and do) assume that $W=W_1$,
$A_1=H_1=0$ (see Assumption I in \S \ref{s-3}), and furthermore that $p=0$. 
By condition (ii) in Theorem \ref{main-r1} (which, as we already pointed out, 
is also assumed in Lemma \ref{spaghetti}) and Lemma \ref{extremepoints} it is 
then enough to prove Lemma \ref{spaghetti} in the case when the index $1$ 
belongs to the set $S(p)$ defined above, which we now proceed to do.

\subsection{Changing Coordinates}\label{ss-coor}

To prove Lemma \ref{spaghetti} in the above situation we will further 
simplify the picture by making suitable coordinate changes as follows. 
Let $G$ be a $C^{1}$-homeomorphism from a domain $U'$
to $U$ that takes a neighborhood $N'\subset U'$ of $p'=G^{-1}(p)$ one-to-one 
onto $N$. Then $W(p)\cap N$ is the homeomorphic image under 
$G$ of the set 
$$W'(p)=\{w\in N'\mid H_{i}(G(w))\leq 0,\,2\le i\le r\}$$
(note that we do not need to assume that $G$ is analytic since 
we are not concerned with preserving subharmonicity in the present situation).
Furthermore, if $z\in U$ and $z'=G^{-1}(z)$ then $V(z)$ is  
the homeomorphic image under $G$ of the set
\begin{equation*}
\begin{split}
&V'(z')\\
&=\{\ze'\in U'\mid \exists\text{ piecewise $C^1$-path from $z'$ to $\ze'$ 
along which all $H_{i}\circ G$ decrease}\}.
\end{split}
\end{equation*} 

Clearly, since $G$ is one-to-one it suffices for the proof of 
Lemma~\ref{spaghetti} to show that there exists a neighborhood $N'$ of 
$p'$ such that 
$V'(z')$ tends to $\mathring{W'}$ through an appropriate set as $z'\to p'$ 
(cf.~Definition \ref{def-lim}). 

As an immediate application of this observation we may prove 
Lemma \ref{spaghetti} in the case when $K$ is a line segment. Indeed, 
suppose that 
$A_1(0)=0$ and $A_2(0)$ are the (only) two extreme points of $K$. By 
Lemma \ref{extremepoints} the functions $A_1(z)\equiv 0$ and $A_2(z)$ 
are the only active ones at $p$ and it suffices to show that $V(z)$ tends to 
$W$ through $W$ as $z\to p$ in a suitable neighborhood.
We may change coordinates as above in order to reduce this 
case to the situation when $H_2(x,y)=y$. Then just consider the harmonic 
conjugate $Q$ of $H_{2}$ and note that $N\ni z\mapsto 
(Q(z),H_{2}(z))$ is a local homeomorphism for a sufficiently small 
neighborhood $N$ of $p=0$. It follows that 
$$V(z)\cap N=\{w\in N\mid \Re w\leq  \Re z\}\text{ and }
\mathring{W}\cap N=\{w\in N\mid \Re w< \Re p=0\},$$ 
so the conclusion of Lemma~\ref{spaghetti} is immediate in this case. 

\subsection{The General Case $r\ge 3$}\label{ss-c1} 

From the discussion at the beginning of this section
it follows that if $W$ is as in~\eqref{w-set} 
and as before $\mathring{W}$ is its interior we
get that the open set 
$$\Om(p):=\bigcap_{i=2}^{r} N_{H_{i}}=\mathring{W}\cap N$$ 
is bounded by 
parts of some of the ($0$-)level curves through $p=0$ of $H_{i}$, 
$2\le i\le r$, and part of the boundary of $N$. Furthermore,
 $\sigma_1(p)$ is a pointed cone subtending 
an angle $\alpha\in (0,\pi)$ at its 
vertex (which is the origin), and it is bounded in a small neighborhood of 
$p$ by tangents to some level curves, say $H_2=0$ and $H_3=0$, that meet 
transversally at $p$. Since two non-identical real 
analytic curves can intersect each other only in a discrete set it follows 
that for a small enough neighborhood $N$ of $p$ the boundary of $\Om(p)$
will consist of at most part of two level curves (and part of the boundary of 
$N$).

By the inverse function theorem the map 
$$(x,y)\mapsto R(x,y):=(H_{2}(x,y),H_{3}(x,y))$$ is a homeomorphism from a 
neighborhood (also called $N$) of $p$ to a 
neighborhood 
of $p$. This map takes $W\cap N $ to an open subset of the 
third quadrant and $p$ is an interior point in the induced topology of the 
third quadrant. Clearly, the homeomorphism $G(x,y)=R^{-1}(x,y)$ satisfies 
$H_{3}(G(x,y))=x$ and $H_{2}(G(x,y))=y$ so that by \S \ref{ss-coor} 
we may assume without loss of generality 
throughout the rest of this section that 
\begin{quote}
{\em $H_{2}(x,y)=y$, $H_{3}(x,y)=x$, $\sigma_1(p)$ is 
the third quadrant and $W\cap N $ is the corresponding quadrant of a disk}.
\end{quote}

The assumption on the boundary of the convex hull of the 
$A_i(p)$'s (cf.~Lemma \ref{extremepoints} and the discussion following it) 
implies that there are no other level curves through 
$p$ that are parallel to either of the level curves of $H_2$ or $H_3$ through 
$p$ except the latter curves themselves. 

Now by viewing gradients as complex numbers for each $z\in N$ we may write 
\begin{equation}\label{grad-com}
\na H_k(z)=\vert \na H_k(z)\vert e^{\sqrt{-1}\theta_k(z)},\text{ where }
\theta_k(z)\in [0,2\pi),\, 2\le k\le r.
\end{equation}
Our assumptions imply 
that $0< \theta_k(p)< \pi/2$ for $2\le k\le r$. Let us further shrink 
$N$ -- if necessary -- so that 
\begin{equation}
\label{eq:theta}
0< \theta_k(z)<\pi/2\text{ for } k\in \{2,\ldots,r\}\setminus \{2,3\} 
\text{ if } z\in N.
\end{equation}

\begin{Claim}\label{new-claim}
For any $z\in \mathring{W}\cap N$ there exists 
a neighborhood $\tilde N_z$ of $0$ such that every point in $\tilde N_z$ may 
be reached by a path from $z$ 
along which each of the functions $H_k$, $2\le k\le r$, increases.
\end{Claim}

\begin{proof}
Let $z\in \mathring{W}\cap N$. Then clearly both coordinates 
$x$ and $y$ are increasing along the straight segment 
from $z$ to $p=0$ given by $\{(1-t)z\mid t\in [0,1]\}$. Moreover, there is 
a disk $N_z$ centered at $p$ such that $w\in N_z$ implies that both 
$x$ and $y$ increase along the path $\gamma_w (t)=(1-t)z+tw$, $t\in [0,1]$, 
from $z$ to $w$. (Note that $N_z$ is the largest disk contained in 
$N\cap \{w\in\bC\mid\Re w\geq \Re z, \Im w\geq \Im z\}$.) Thus both functions
$[0,1]\ni t\mapsto H_k(\gamma_w(t))$, $k\in\{2,3\}$, are increasing. 
Let us show that this is true as well for each of the remaining functions 
$[0,1]\ni t\mapsto H_k(\gamma_w(t))$, $k\in \{2,\ldots,r\}\setminus \{2,3\}$.
By~\eqref{eq:theta} one has $\na H_k(z)=(\alpha(z),\beta(z))$, 
where $\alpha(z),\beta(z)>0$ if $k\notin\{2,3\}$ and $z\in N$, so that
the derivative 
\begin{equation}
\label{eq:riktnder.lin}
\frac{d}{dt}H_k(\gamma_w(t))
=\alpha(\gamma_w(t))\Re(w-z)+\beta(\gamma_w(t))\Im(w-z)
\end{equation}
is positive for $w=0$, $2\le k\le r$, and $t\in[0,1]$. Hence there is 
a neighborhood $\tilde N_z$ of $0$ such that the expression 
in~\eqref{eq:riktnder.lin} is positive for all $w\in \tilde N_z$ and 
$t\in[0,1]$. 
This means that each point in $\tilde N_z$ may be reached by a path from $z$ 
along which each of the functions $H_k$, $2\le k\le r$, increases.
\end{proof} 

The proof of Lemma~\ref{spaghetti} is now immediate: 
if $\{z_n\}_{n\in \bN}$ is a sequence converging to $p$ there is 
$n_0\in \bN$ such that $n\geq n_0$ implies $z_n\in \tilde N_z$ and by 
Claim \ref{new-claim} 
there is a path from $z$ to $z_n$ along which all $H_k$, $2\le k\le r$, 
increase. Going in the other direction there is a path from $z_n$ to $z$ 
along which all $H_k$, $2\le k\le r$, decrease hence $z\in V(z_n)$ for 
$n\geq n_0$. By the above remarks this 
completes the proof of Lemma~\ref{spaghetti}.

\subsection{A More Precise Version of Theorem \ref{main-r1}}\label{ss-44}

Revisiting the proof of Theorem~\ref{main-r1} sketched in \S \ref{s-3} 
we see that we can 
actually formulate a more precise result by using the terminology and 
arguments given in \S \ref{ss-prel}--\S \ref{ss-c1} above.

\begin{Cor}\label{Cor5}
Assume that all hypotheses of Theorem \ref{main-r1} are satisfied {\em except  
condition (ii)}. Let $A_i(p)$ be an extremal point in $K$ 
and consider the part $\partial K_i$ of its boundary (i.e., the union of the 
two edges of $K$) connecting $A_i(p)$ to its two neighbouring extremal points. 
If $A_k(p)\notin \partial K_i$, $k\neq i$, 
there exists a neighborhood $N$ of $p$ such that 
$\Psi=2\partial \varphi /\partial z$ a.e.~in $W_{i}(p)\cap N$.
\end{Cor}
    
\section{Proof of Lemma~\ref{lemma:decreasing}}\label{s-5}

In this section we prove the remaining lemma, namely 
Lemma~\ref{lemma:decreasing} that generalizes a 
corresponding result obtained in  \cite{BR} in the (simpler) case when the 
$A_i$ are constant functions. 
Recall Notation~\ref{not-2}, the renormalization argument in Assumption I of
\S \ref{s-3} allowing $A_1\equiv 0$, and the assumptions of 
Lemma~\ref{lemma:decreasing} and Theorem \ref{main-r1}
for our given $PA$ function 
\begin{equation}\label{eq-norm}
\Psi(z)=\sum_{i=1}^{r}A_i(z)\chi_i(z)=0\cdot \chi_1(z)
+\sum_{i=2}^{r}A_i(z)\chi_i(z)
\end{equation}
and for the path $\gamma$. In particular, we assume that condition (iii) in 
Theorem \ref{main-r1} is fulfilled at all points on $\gamma$, that is,
$\gamma$ does not pass through singular points for the differences $H_i-H_j$
with $i\neq j$.  
We may reparametrize $\gamma$ by arc-length using the parameter interval 
$[0,L]$ and so we may assume that $\vert {\dot{\gamma}(t)}\vert =1$, 
$t\in [0,L]$. Note first that it is enough to prove the following modified 
form of Lemma~\ref{lemma:decreasing}: for each $t_1\in [0,L]$ there exists 
$\eta>0$ such that for any positive test function $\phi$ with 
$\text{supp}\,\phi$ small enough one has
\begin{equation}\label{chi:test2}
\begin{split}
&(\chi_{1}*\phi)(z_{1})\leq (\chi_{1}*\phi)(z_{2}),\\
&\text{where } z_1=\gamma(t_1) \text{ and } z_2=\gamma(t_2)
\text{ with } 0<t_2-t_1<\eta.
\end{split}
\end{equation}
Indeed, the fact that~\eqref{chi:test2} implies Lemma~\ref{lemma:decreasing} 
follows easily by a compactness argument: fix $t_1$ and let $s_2$ be maximal 
such that \eqref{chi:test} holds for $t_2<s_2$. If $s_2\neq L$ 
then~\eqref{chi:test2} gives a contradiction to the maximality of $s_2$. 
For simplicity we make a translation so that $z_1=0$. Clearly, we may also 
assume that $\gamma$ is $C^1$.

The idea of the proof of inequality~\eqref{chi:test2}
is to use the asymptotic properties of the logarithm of $\Psi$. 
For this we need to take the logarithm of the $A_i$ and we must 
therefore make sure that it is possible to choose a suitable branch.
To this end we first prove the following assertion.

\begin{Claim}\label{cl-bla}
There exists a neighborhood $M$ of $z_1=0$ such that 
$$A_{i}(z)\in\bC\setminus \{t \bar{v}\mid t\in (0,\infty)\},\quad 
z\in M,\,1\le i\le r,$$
whenever $v$ is a unimodular complex number satisfying $v\in\sigma(z_1)$, 
where (cf.~\eqref{eq-cone})
$$\sigma(z_1)=\bigcap_{i=2}^{r} \{u\mid \Re \left[u A_i(z_1)\right]\leq 0\}.$$
\end{Claim}

\begin{proof}
Since $A_1\equiv 0$ this is immediate for $i=1$.
By condition (iii)  in Theorem~\ref{main-r1} there exists
$c'>0$ such that $|A_{i}(z_1)|\ge c'$ for $i\in \{2,\ldots,r\}$, 
so that there is
$c\in (0,c']$ and a neighborhood $M$ of $z_1$ such that $|A_{i}(z)|\ge c$ for  
$i\in \{2,\ldots,r\}$ and $z\in M$. 
It follows that for all unit vectors 
$v\in \sigma(z_1)$ we may assume up to shrinking $M$ that 
$\Re \left[vA_{i}(z)\right]\leq c/2$ for $z\in M$. Thus the angle $\rho$
between $A_{i}(z)$ and $\bar v$ satisfies $\rho\in(\pi/3, 5\pi/3)$ since 
$\cos\rho =|A_{i}(z)|^{-1}\Re\left[vA_{i}(z)\right]<1/2$,  
which proves the claim.
\end{proof}

We use the result that we have just established in order 
to simplify the situation. For this we choose $\eta>0$ such that
$\gamma(t)\in M$, $t\in [0,\eta]$, where the neighborhood $M$ of $z_1=0$ is 
as in Claim \ref{cl-bla}, and we 
let $v= {\dot{\gamma}(0)}$. Note that since by the assumption in 
Lemma~\ref{lemma:decreasing} all functions 
$[0,\eta]\ni t\mapsto H_i(\gamma(t))$,
$2\le i\le r$, are decreasing we have $v\in \sigma(z_1)$ 
by \eqref{eq:riktnder}.
Up to replacing $\Psi$ by the function
$e^{i\theta}\Psi(e^{i\theta}z)$,
where $v=e^{i\theta}$, we may also assume that $v=1$. 
In particular, we deduce that $\Re\left[{\dot{\gamma}(0)}\right]=1>0$ so 
that by further shrinking $M$ and the corresponding $\eta>0$ 
we get the key property
\begin{equation}\label{Re}
\Re\left[{\dot{\gamma}(t)}\right]>0,\quad t\in [0,\eta].
\end{equation}
 
Let $\widetilde \Psi_{\epsilon}=\log(\Psi-\epsilon)$,
where $\epsilon>0$ is
arbitrary and we have chosen a branch of the logarithm that is
defined in the complex plane cut along the positive real axis.
The composite distribution $\tilde \Psi_{\epsilon}$ is then 
defined by the above rotation of the complex plane, since 
$v=1\in \sigma(z_1)$. We now study its derivative along the path $\gamma$. 

Give $\ze\in M$ define as above (cf.~\eqref{eq-cone}) 
$$\sigma(\ze)=\bigcap_{i=2}^{r} \{u\mid \Re \left[u A_i(\ze)\right]\leq 0\}.$$
Then for any fixed $\epsilon >0$ and $u\in \sigma(\ze)$ with $\Re u>0$
one has 
\begin{equation}\label{villkor}
\Re\left[u(A_{i}(w)-{\epsilon})\right]<0,\quad 1\le i\le r,
\end{equation}
for all 
$w$ in a (sufficiently small) neighborhood of $\ze$. In particular, 
inequality \eqref{villkor} holds 
for all vectors of the form $u={\dot{\gamma}(t)}$ in view of \eqref{Re} 
and the fact that all functions $[0,\eta]\ni t\mapsto H_i(\gamma(t))$,
$2\le i\le r$, are decreasing (and thus $u\in \sigma(\ze)$ 
by \eqref{eq:riktnder}). It follows that if $\phi$
is a positive test function with $\int\!\phi d\lambda=1$ and 
$\text{supp}\,\phi$ is small enough then
\begin{equation}\label{villkor2}
\Re\left[u(\phi*\Psi-\epsilon)\right]< 0
\end{equation}
and therefore 
$$\Re\left[\frac{\bar u}{\phi* \Psi-\epsilon}\right]\leq 0$$
in a neighborhood of $\ze$. 
Since $\partial(\phi*\Psi)/\partial \bar z\geq 0$ we get 
$$\Re\left[\bar u\frac {\partial} {\partial {\bar z}}
{\log (\phi* \Psi-\epsilon)}\right]=
\Re\left[\frac{\bar u}{\phi* \Psi-\epsilon }\cdot
\frac{\partial(\phi* \Psi)} {\partial {\bar z}}\right]\leq 0.$$
Letting $\text{supp}\,\phi\to 0$ with
$\int\!\phi d\lambda=1$ we see that 
$\log (\phi* \Psi-\epsilon)\to \widetilde \Psi_{\epsilon}$ in $L^1_{loc}$
(hence as a distribution) and by passing to the limit we get
$$\Re\left[\bar u\frac{\partial\widetilde \Psi_{\epsilon}} 
{\partial{\bar z}}\right]\leq 0.$$
Write now $\widetilde \Psi_{\epsilon}=\sigma_{\epsilon}+i\tau_{\epsilon}$, 
where $\sigma_{\epsilon}$ and $\tau_{\epsilon}$ are real-valued distributions.
Then the latter inequality yields
\begin{equation}\label{Reimlog}
\Re\left[\bar u\frac{\partial\sigma_{\epsilon}}{\partial \bar z}\right]\leq
\Im\left[\bar u\frac{\partial\tau_{\epsilon}}{\partial\bar z}\right],
\end{equation}
where \eqref{Reimlog} is interpreted as being valid for the restrictions of 
the corresponding distributions to a neighborhood of $\ze$. Note that
up to further shrinking $M$ (and the corresponding $\eta>0$) 
by our choice of the branch of the logarithm used in the definition of
$\widetilde \Psi_{\epsilon}$ we have
\begin{equation}\label{b-im}
\tau_\eps(z)\in \left(\frac{\pi}{2},\frac{3\pi}{2}\right),\quad z\in 
2M=\{a+b\mid a,b\in M\}.
\end{equation}

Let us show that relations \eqref{Reimlog}--\eqref{b-im} 
produce the desired result. Recall that 
for a real-valued function $\omega(z)$ one has 
\begin{equation}\label{bla-om}
\frac{\partial \omega(z)}{\partial\bar z} =\overline{\frac{\partial 
\omega(z)}{\partial z}}
\end{equation}
in the sense of distributions. We consider the derivative of 
$\widetilde \Psi_{\eps}$ along the path $\gamma$: if $\phi$
is a positive test function  
then since $\sigma_{\epsilon}$ is a real-valued distribution we deduce
from \eqref{bla-om} and \eqref{Reimlog} that 
the following holds in the interval $(0,\eta)$:
\begin{equation}\label{bla1}
\begin{split}
\frac{d}{dt}\left[(\phi* \sigma_{\epsilon})(\gamma (t))\right]
&=2\Re\left[\dot{\gamma}(t)\frac{\partial\phi*
\sigma_{\epsilon}}{\partial z} (\gamma(t))\right]
=2\Re\left[\overline{\dot{\gamma}(t)}
\frac{\partial\phi*\sigma_{\epsilon}}{\partial\bar z}
(\gamma(t))\right]\\
&=2\int\Re\left[\overline{\dot{\gamma}(t)}
\frac{\partial\phi}{\partial\bar z}(\gamma(t)-w)\sigma_{\epsilon}(w)
\right]d\lambda(w)\\
&\leq 2\int\Im\left[\overline{\dot{\gamma}(t)}
\frac{\partial\phi}{\partial\bar z}(\gamma(t)-w)\tau_{\epsilon}(w)
\right]d\lambda(w).
\end{split}
\end{equation}
Now if $\text{supp}\,\phi$ is small enough, say $\text{supp}\,\phi
\subset M$,
then from \eqref{b-im} and the fact that $|\dot{\gamma}(t)|=1$ for 
$t\in [0,\eta]$
(cf.~the reparametrization argument at the beginning of this section) we get
\begin{equation}\label{bla2}
2\left|\int\Im\left[\overline{\dot{\gamma}(t)}
\frac{\partial\phi}{\partial\bar z}(\gamma(t)-w)\tau_{\epsilon}(w)
\right]d\lambda(w)\right|\le 2\cdot \frac{3\pi}{2}\cdot\frac{1}{2}
\left(\left|\!\left|\frac{\partial\phi}{\partial x}\right|\!\right|_1+
\left|\!\left|\frac{\partial\phi}{\partial y}\right|\!\right|_1\right)=:
\kappa(\phi),
\end{equation}
where $|\!|\cdot|\!|_1$ denotes the $L^1$-norm. Note that the (positive)
constant $\kappa(\phi)$ defined above does not depend on $\eps$. 
Combining \eqref{bla1} and \eqref{bla2} we obtain
\begin{equation}\label{9}
(\phi* \sigma_{\epsilon})(z_{2})-
(\phi*\sigma_{\epsilon})(z_{1})\leq \kappa(\phi)\eta.
\end{equation}
On the other hand by \eqref{eq-norm} we have
$$\tilde{\Psi}_\eps(z)
=\log\!\left[-\eps\chi_1(z)+\sum_{i=2}^{r}(A_i(z)-\eps)\chi_i(z)\right]$$
hence
$$\sigma_\eps(z)=(\log\eps)\cdot\chi_1(z)+f_\eps(z),\,\text{ where }\,
f_\eps(z)=\sum_{i=2}^{r}\log|A_i(z)-\eps|\cdot
\chi_i(z),$$
and therefore $(\phi* \sigma_{\epsilon})(z)=(\log\epsilon)\cdot
(\phi*\chi_{1})(z)+(\phi* f_\eps)(z)$.  
By condition (iii)  in Theorem \ref{main-r1} there exists $c>0$ 
such that $|A_{i}(z)|\ge c$ for  
$i\in \{2,\ldots,r\}$ and $z\in M$ (cf.~the proof of Claim \ref{cl-bla}). We 
deduce that there exists $c'>0$ (independent of $\eps$ and $\phi$) such that 
$|(\phi* f_\eps)(z)|\le c'|\!|\phi|\!|_\infty$ for $z\in M$, where 
$|\!|\cdot|\!|_\infty$ denotes the $L^\infty$-norm. It follows that 
\begin{equation}\label{10}
%\begin{split}
(\phi* \sigma_{\epsilon})(z)
=(\log\epsilon)\cdot(\phi*\chi_{1})(z)
%\sum_{i=2}^{r}\log|A_i(z)-\eps|\cdot(\phi*\chi_i)(z)\\
%&=(\log\epsilon)\cdot(\phi*\chi_{1})(z)
+O(1).
%\end{split}
\end{equation}
Substituting \eqref{10} in \eqref{9} and letting $\epsilon \to 0$ we conclude 
that~\eqref{chi:test2} holds, which by the preliminary remarks 
at the beginning of this section completely settles 
Lemma~\ref{lemma:decreasing}.

\section{An Alternative Approach Under Extra Conditions}\label{s-ref}

In the previous sections we formulated and proved three results answering 
the Main Problem stated in \S \ref{s-1} under fairly mild assumptions, 
namely Theorem \ref{main-r1} and its consequences Corollary \ref{Cor5} and 
Corollary \ref{cor-ip-conv} (cf.~\eqref{eq-chain}). We will now prove
Theorem~\ref{th-ref} below that provides a fourth answer to the Main Problem 
under some extra (yet still mild) conditions. Although this result 
may be obtained directly from Corollary~\ref{cor-ip-conv}, the point in 
what follows is to present a different approach\footnote{This approach and the 
subsequent proofs were suggested by the referee whom we would like to 
thank for generously sharing his ideas with us.} from the one used in 
\S \ref{s-3}--\S \ref{s-5} that does not rely on Lemma \ref{lemma:decreasing} 
and Lemma~\ref{spaghetti}. 

\begin{notation}\label{not-ref}
Let $\Phi\in PA$ be as in \eqref{eq:MDA}, which we assume to be the canonical 
piecewise decomposition of $\Phi$ in the sense of Definition \ref{d-r2}. 
We may write
$$
U\setminus\bigcup_{i=1}^{r}M_i=Z,
$$
where $M_i$, $1\le i\le r$, are pairwise disjoint open sets and $Z$ 
is Lebesgue negligible. Note that each $\pa M_i$ is also Lebesgue negligible 
since $\pa M_i\subset Z$, $1\le i\le r$. As before we let $\chi_i$ be the 
characteristic 
function of $M_i$. Recall from \eqref{ip} the set $I(p)$ and its cardinality 
$i(p)$ defined for any $p\in U$. To simplify some discussions, assume that 
$U$ is simply connected and 
choose $f_i\in A(U)$ such that $f_i'(z)=A_i(z)$, $1\le i\le r$, where the 
$A_i$ are the given (analytic) functions appearing in the decomposition 
\eqref{eq:MDA} of $\Phi$. Hence
$$
\Phi(z)=\sum_{i=1}^{r}f_i'(z)\chi_i(z).
$$
For arbitrarily fixed $p\in U$ we let 
$$
\phi(z)\,(=\phi_p(z))
=\max_{j\in I(p)}\Re(f_j(z)-f_j(p))=\max_{j\in I(p)}H_j(z),
$$
where the $H_i$ are the harmonic functions defined in \eqref{a-h} 
(cf.~\S \ref{s-1} in the case when $i(p)=r$). 
Clearly, $\phi$ is a continuous subharmonic function in $U$ which vanishes 
at $p$. Finally, if $k\in I(p)$ and $i(p)>1$ set
$$
V_k(p)=\left\{\sum_{j\in I(p)\setminus \{k\}}
\te_j\left(f_k'(p)-f_j'(p)\right)\,\Bigg|\, 
\te_j\ge 0,\,j\in I(p)\setminus \{k\},\,\sum_{j\in I(p)\setminus \{k\}}
\te_j>0\right\}.
$$
\end{notation}

Recall the definition of $\Si(U)$ from Notation \ref{not-Sigma}.

\begin{Theorem}\label{th-ref}
In the above notations assume that $\Phi\in\Si(U)$ and that the following 
conditions hold:
\begin{itemize}
%\item[(i)] Each $\pa M_i$ is of Lebesgue measure $0$;
\item[(i)] The one-dimensional Hausdorff measure of 
$\pa M_j\cap \pa M_k \cap \pa M_l$ is $0$ whenever $j<k<l$;
\item[(ii)] If $i(p)>1$ and $k\in I(p)$ then $0\notin V_k(p)$.
\end{itemize}
Then a.e.~in a neighborhood of every $p\in U$ one has 
$\Phi=2\pa_z\phi\,(=2\pa_z\phi_p)$.
\end{Theorem}

\begin{remark}\label{rem-ref}
Recall the assumptions involving the (extremal points of the) convex hull $K$ 
of the points $A_i(p)=f_i'(p)$, $1\le i\le i(p)$, that we used in 
Corollary \ref{cor-ip-conv} and Corollary \ref{Cor5}. Although still mild 
(since it is generically true), requirement (ii) in Theorem \ref{th-ref} is 
actually stronger than the aforementioned assumptions.
\end{remark}

The remainder of this section is devoted to the proof of Theorem \ref{th-ref},
which uses induction on $i(p)$. 

Consider first the case $i(p)=1$. By relabeling the indices we may assume 
that $I(p)=\{1\}$, that is, $p\notin \overline{M}_j$ for $j>1$. Hence $p$ is 
either an interior point of $M_1$ or $p\in Z$ and every neighborhood of 
$p$ intersects $M_1$. If the former occurs then
$\Phi(z)=2\pa_z\Re f_1(z)$ in an open neighborhood of $p$ and thus 
$\Phi=2\pa\phi/\pa z$
in that neighborhood. If $p\in Z$ then there is a small open neighborhood 
$\Om$ of $p$ contained in $M_1\cup Z$ and we conclude that 
$\Phi(z)=2\pa_z\Re f_1(z)$ a.e.~in $\Om$, hence equality holds in $\Om$ 
in the distribution sense. This settles the case when $i(p)=1$.

Assume next that $i(p)=2$ and (without loss of 
generality) $I(p)=\{1,2\}$. Since the $M_i$ are pairwise disjoint it follows 
that $p\in Z$ and $p\notin \overline{M}_k$ for $k>2$. Therefore, there is 
an open neighborhood $\Om$ of $p$ such that
$$
\Phi(z)=f_1'(z)\chi_1(z)+f_2'(z)\chi_2(z),\quad z\in\Om.
$$
Let $\chi=\chi_2\big|_\Om$, $f=f_2-f_1\big|_\Om$, and define
$$
\Psi(z)=f'(z)\chi(z)=\Phi(z)-f_1'(z).
$$
Note that $\pa_{\bar{z}}\Psi(z)\ge 0$ in $\Om$. Condition (ii)
in Theorem \ref{th-ref} implies that $f'(p)\neq 0$ and we may assume 
(after shrinking $\Om$, if necessary) that $f$ is a diffeomorphism from $\Om$ 
onto some open disk $D\subset \bC$. We may then write $\chi(z)=\eta(f(z))$, 
where $\eta=\eta(w)=\eta(u+iv)$ is the characteristic function of some 
open subset $\omega$ of $D$, and we get
$$
0\le \pa_{\bar{z}}\Phi(z)=\pa_{\bar{z}}f'(z)\eta(f(z))
=|f'(z)|^2(\pa_{\bar{w}}\eta)(f(z)),
$$
so that $\pa_{\bar{w}}\eta\ge 0$ in $D$. Since $\eta$ is real-valued this 
means that $\eta$ is an increasing function of $u$. Hence 
the open set 
$\omega$ is defined by an inequality of the form $\Re w>a$, and then 
$M_2\cap \Om$ is defined by the inequality $\Re(f_2(z)-f_1(z))>a$. Moreover, 
since $p$ is in the closure of the set where $\chi=1$ we must have 
$a=\Re(f_2(p)-f_1(p))$. Clearly, we may assume that $f_1(p)=f_2(p)=0$. 
Then $\Phi(z)=f_1'(z)$ when $z\in\Om$ and $\Re f_1(z)>\Re f_2(z)$ while 
$\Phi(z)=f_2'(z)$ when $z\in\Om$ and $\Re f_1(z)<\Re f_2(z)$. This 
shows that $\Phi=2\pa\phi/\pa z$ in a neighborhood of $p$, which completes
the proof in the case when $i(p)=2$.

The above observations also give us a result that will be used 
later on:

\begin{Lemma}\label{ref-l1}
Assume that $I(p)=\{j,k\}$, where $j<k$, and that $\ga(t)$ is a $C^1$-curve 
escaping from $M_j$ into $M_k$ when $t=\tau$ in the sense that $\ga(t)\in M_j$ 
for $t<\tau$ and there is a sequence $\{\tau_\nu\}_{1}^{\infty}$ with 
$\tau_\nu>\tau$ and $\tau_\nu\to \tau$ as $\nu\to \infty$ 
such that $\ga(\tau_\nu)\in M_k$. 
Then $\pa_t\Re\big(f_j(\ga(t))-f_k(\ga(t))\big)\Big|_{t=\tau}\le 0$.
\end{Lemma}

Let us now pass to the case when $i(p)\ge 3$. Then $p\in Z$ and there is 
an open neighborhood of $p$ that does not intersect $r-i(p)$ of the 
$\overline{M}_j$. By deleting these sets from $U$ we may assume that 
$i(p)=r\ge 3$ (cf.~the comments after \eqref{ip} in \S \ref{s-1}). We then know
that $p\in\bigcap_{i=1}^{r}\pa M_j$. It is no restriction to further assume 
that the $f_j$ are normalized so that $f_j(p)=0$ for every $j$. Then 
$\phi(z)\,(=\phi_p(z))=\max_j\Re f_j(z)$ and we have to prove that
\begin{equation}\label{eq-71}
\Re f_k=\phi\,\text{ in }\,M_k\cap N,
\end{equation}
where $N\subset U$ is a sufficiently small open neighborhood of $p$. Let
$$
N_k=\{z\in N\mid \Re f_k(z)>\Re f_j(z)\text{ when }j\neq k\}.
$$

Suppose now that we can show the following:
\begin{equation}\label{eq-72}
N_k\subset \overline{M}_k\,\text{ for every $k$ if $N$ is sufficiently small}.
\end{equation}
Since the $\Re f_j$ must be pairwise distinct harmonic functions in $U$ (as 
a consequence of condition (ii) in Theorem \ref{th-ref}), the set where 
$\Re f_j=\Re f_k$ for some $j,k$ with $j\neq k$ is of Lebesgue measure $0$. 
It follows that $N$ is the disjoint union of the sets $N_k$ together with a 
set of measure $0$. Since the $M_j$ are pairwise disjoint and $\pa M_j$ is of 
Lebesgue measure $0$ for every $j$ (since $\pa M_j\subset Z$, 
cf.~Notation \ref{not-ref}) we deduce that 
$(M_k\cap N)\setminus N_k$ is Lebesgue negligible. From this we conclude that 
$\Re f_k=\phi$ in $M_k\cap N$ hence $\Phi=2\pa_{z}\phi$ in $N$, which proves 
Theorem \ref{th-ref}.

Thus the main issue is to show that \eqref{eq-72} holds. 
When doing this we may assume 
that $k=r$ and consider the harmonic functions $h_j=\Re(f_r-f_j)$, 
$1\le j\le r-1$. We know that $h_j(p)=0$. Let $q\in N_r$, i.e., $q\in N$ and 
$h_j(q)>0$ for $j<r$. We want to show that $q\in \overline{M}_r$. For this 
we define
$$
\La=\bigcup_{j<k<l}\left(\pa M_j\cap\pa M_k\cap \pa M_l\right).
$$
By assumption (i) in Theorem \ref{th-ref}, $\La$ has vanishing one-dimensional
Hausdorff measure.
We need the following lemma.

\begin{Lemma}\label{ref-l2}
There is an open set $N\subset U$ containing $p$ such that the following 
holds: if $w\in N$ and
$h_k(w):=\Re(f_r(w)-f_k(w))>0$ when $k<r$,
then there exist an open neighborhood $\cM=\cM_w\subset U$ of $p$ and for 
every $z\in \cM$ a real analytic mapping $\ga=\ga(s,t)$ from a neighborhood of 
$[0,1]\times [0,1]$ into $U$ such that
\begin{itemize}
\item[(a)] The restriction of $\ga$ to any set where $t<t_0<1$ is a 
diffeomorphism onto its image;
\item[(b)] $\ga(1/2,0)=z$ and $\ga(s,1)=w$ for all $s$;
\item[(c)] $\pa_t h_k(\ga(s,t))>0$ for all $(s,t)$ 
when $k<r$.
\end{itemize}
\end{Lemma}

Assertion \eqref{eq-72} -- and thus, as explained above, Theorem \ref{th-ref} 
as well -- is now a 
consequence of Lemma \ref{ref-l2}. Indeed, let $N$ be a small neighborhood 
of $p$ satisfying the assumptions of Lemma \ref{ref-l2} and $w\in N$ be 
such that $h_k(w)>0$ for $k<r$. We need to prove that $p\in\overline{M}_r$. 
For this let $\cM=\cM_w$ be as in the conclusion of Lemma \ref{ref-l2}. 
Since $p\in\overline{M}_r$ we know that $\cM$ contains a point $z\in M_r$. 
Let $\ga$ be the mapping corresponding to $z$ and $w$. By shrinking the domain 
in which the variable $s$ ranges we may assume that $\ga(s,0)\in M_r$ when 
$s\in [0,1]$. Set
$$
\cA_\nu=\left\{(s,t)\mid0\le s\le 1,\,0\le t\le 1-\nu^{-1}\right\}
$$
for each integer $\nu\ge 2$. Since the one-dimensional
Hausdorff measure of $\La$ vanishes this is also true for the one-dimensional
Hausdorff measure of 
$$
K_\nu:=\{(s,t)\in\cA_\nu\mid \ga(s,t)\in\La\}.
$$
It follows that
$$
J_\nu:=\{s\in[0,1]\mid (s,t)\in K_\nu\text{ for some }t\}
$$
is a closed set of Lebesgue measure $0$. In fact, $J_\nu$ is the projection 
of a set with vanishing one-dimensional
Hausdorff measure, see, e.g., \cite[Theorem 7.5]{Mat}. Therefore, the set 
$J_\nu$ is of the 
first category, which implies that $\bigcup_\nu J_\nu$ is also of the first 
category. This gives us an $s\in[0,1]$ such that $\ga(s,t)\notin \La$ 
when $0\le t<1$. From condition (c) in Lemma \ref{ref-l2} and 
Lemma \ref{ref-l1} 
it follows that the curve $t\mapsto \ga(s,t)$, which starts at 
$\ga(s,0)\in M_r$, can not leave $\overline{M}_r$ until $t=1$. Hence 
$w\in \overline{M}_r$, which proves \eqref{eq-72} and we are done.

It remains to prove Lemma \ref{ref-l2}. In doing so we will use the fact that 
the functions $h_j=\Re(f_r-f_j)$, $1\le j\le r-1$, introduced above are 
real-valued
and real analytic, but we will make no use of their harmonicity. Condition 
(ii) in Theorem \ref{th-ref} implies that the set of all linear combinations 
$\sum_{j=1}^{r-1}\te_jdh_j(p)$, where $\te_j\ge 0$ for all $j$ and $dh$ 
denotes differential, is contained in a convex cone $\Ga$ with positive 
opening angle less than $\pi$. We make an affine change of coordinates, 
only keeping the affine space structure of $\bC$. This change of coordinates 
will allow us to replace $\Ga$ with any other cone with positive opening angle,
and without loss of generality we may further assume that $p$ is the origin.
Then we are in the situation where a set of $m=r-1$ real analytic and 
real-valued functions $h_1,\ldots,h_m$ are defined in a neighborhood $V$ of 
the origin in $\bR^2$ and satisfy the conditions
\begin{itemize}
\item[(I)] $h_j(0)=0$ and $dh_j(0)\neq 0$ when $1\le j\le m$; 
\item[(II)] The closed convex cone generated by the gradients $\nabla h_j(0)$, 
$1\le j\le m$, is contained in the cone
$\Ga:=\{(x,y)\in\bR^2\mid |x|\le y\}.$
\end{itemize}

To complete the proof of Lemma \ref{ref-l2} we only have to 
establish the following result.

\begin{Lemma}\label{ref-l3}
Assume conditions (I)--(II) above. Then there is an open set $0\in N\subset V$ 
such that the following holds: if
$$
w\in\Om_N:=\{z=(x,y)\in N\mid h_j(z)>0,\,1\le j\le m\}
$$
one can find an open neighborhood $\cM=\cM_w$ of the origin and for each 
$z\in\cM$ a $C^1$-mapping $\ga(s,t)$ from a neighborhood of $[0,1]\times [0,1]$
into $V$ such that
\begin{itemize}
\item[(a)] The restriction of $\ga$ to any set where $t<t_0<1$ is a 
diffeomorphism onto its image;
\item[(b)] $\ga(1/2,0)=z$ and $\ga(s,1)=w$ for all $s$;
\item[(c)] $\pa_t h_k(\ga(s,t))>0$ for all $(s,t)$ 
when $k\le m$.
\end{itemize}
\end{Lemma}

\begin{proof}
Define
$$
\Om_{N}^{\pm}=\Om_N\cap \{(x,y)\in\bR^2\mid \pm x\ge 0\}
$$
whenever $N\subset V$. It suffices to prove that there exist an open set  
$0\in N=N_+\subset V$ such that the conclusion of the lemma holds when 
$w\in\Om_N^+$. Indeed, by replacing $h_k(x,y)$ with $h_k(-x,y)$ we would 
obtain $N=N_{-}$ for which the conclusion of the lemma would then be true when 
$w\in\Om_N^-$ and thus the assertions in the lemma would follow for 
the open set $N=N_+\cap N_{-}$.

It is no restriction to assume that $dh_j(0)$ is proportional to 
$-dx+dy$ for some $j$. By shrinking $V$ if necessary and applying the implicit 
function theorem we may also assume that every $h_j$ is of the form 
$$
h_j(x,y)=\be_j(x,y)(y-g_j(x)),
$$
where $\be_j,g_j$ are real analytic functions and $\be_j>0$. Then by using 
the real analyticity of the functions $g_j$ we may further assume -- after 
shrinking $V$ and relabeling the indices, if necessary -- that 
$V=(-b,b)\times (-b,b)$ for some positive real number $b$ and that 
$g_1(x)\le g_2(x)\le\cdots\le g_m(x)$ when $0<x<b$. With these normalizations 
it follows that
$$
-1\le g_1'(0)\le g_2'(0)\le\cdots \le g_m'(0)=1
$$
and finally, after making a non-linear change of the $x$-coordinate, we may 
additionally assume that $g_m(x)=x$.

Below we let $a<b$ and $\de$ be small positive numbers and we make generic 
use of the letter $C$ to denote constants that are independent of $a$ and 
$\de$ when these stay small. Define
\begin{equation*}
\begin{split}
N(a)&=\{z\in\bC\mid |z|<a\},\\
\Om^+(a)&=\{z=(x,y)\in N(a)\mid x\ge 0\text{ and }h_k(z)>0\text{ for all }k\}, 
\end{split}
\end{equation*}
so that 
$\Om^+(a)=\{z=(x,y)\mid 0\le x <y,\,|z|<a\}.$

Now, we clearly have the estimates
\begin{equation}\label{eq-73}
C^{-1}\le \be_j(z)\,\text{ and }\, |\nabla\be_j(z)|\le C,\quad z\in N(a).
\end{equation}
Let $w=(u,v)\in\Om^+(a)$ and set $\rho=v-u$. Then $\rho$ is a positive real 
number that depends on $w$ and we define
$$
\cM=\cM_w=\{z\in\bC\mid |z|<\de\rho\}.
$$
Take $z\in\cM$ and let $\al\in\bR^2$ be linearly independent from $w-z$ and 
such that $|\al|\le \de\rho$. Introduce the mapping
\begin{equation}\label{extra1}
\ga(s,t)=(x(s,t),y(s,t))=z+(s-1/2)(1-t)\al+t(w-z)
\end{equation}
defined for all $(s,t)$ in a small open neighborhood of $[0,1]\times[0,1]$. It 
is then immediate that assertions (a) and (b) in the lemma are satisfied.

In order to verify (c) we compute the $t$-derivative of $h_j(\ga(s,t))$:
\begin{equation}\label{eq-74}
\begin{split}
\pa_t(h_j(\ga(s,t)))=(y(s,t)-&g_j(x(s,t)))\pa_t(\be_j(\ga(s,t)))\\
&+\be_j(\ga(s,t))(\pa_t y(s,t)-g_j'(x(s,t))\pa_t x(s,t)).
\end{split}
\end{equation}
We see that
\begin{equation}\label{eq-75}
\left|\pa_t(\be_j(\ga(s,t)))\right|\le Ca.
\end{equation}
Since $g_j(x)\le g_m(x)=x$ when $0<x<a$ we may write
$$
g_j(x)=x-p_j(x),
$$
where $p_j(x)\ge 0$. If $p_j(x)\not\equiv 0$ then $p_j(x)=x^{\mu_j}q_j(x)$, 
where $\mu_j$ is a positive integer and $q_j(0)>0$. By taking $a$ sufficiently 
small we may then assume that
\begin{equation}\label{eq-76}
p_j'(x)=\mu_jx^{\mu_j-1}q_j(x)+x^{\mu_j}q_j'(x)\ge C^{-1}p_j(x)/x,\quad 
0<x<a.
\end{equation}
Moreover, since $x(s,t)=(1-t)x(s,0)+tx(s,1)\ge (1-t)x(s,0)$ it follows that 
$|x(s,t)|\le C\de\rho$ if $x(s,t)\le 0$. Hence there is a constant $C$ such 
that
\begin{equation}\label{eq-77}
\left|p_j'(x(s,t))-p_j'(|x(s,t)|)\right|\le C\de\rho,\quad 0\le s,t\le 1.
\end{equation}
Next, one has 
\begin{equation}\label{extra2}
\begin{split}
y(s,t)-&g_j(x(s,t))=(1-t)y(s,0)+ty(s,1)-x(s,t)+p_j(x(s,t))\\
&=(1-t)y(s,0)+ty(s,1)-(1-t)x(s,0)-tx(s,1)+p_j(x(s,t))\\
&=(1-t)(y(s,0)-x(s,0))+t(y(s,1)-x(s,1))+p_j(x(s,t))\\
&=(1-t)(y(s,0)-x(s,0))+t\rho+p_j(x(s,t)).
\end{split}
\end{equation}
Recall that $w\in\Omega^+(a)$, so that in particular $|w|<a$. Since 
$|z|<\de\rho$ and $|\al|\le \de\rho$ it follows from \eqref{extra1} that 
$|x(s,t)|<a$ if $\de$ is small enough. We then deduce from \eqref{eq-77} 
and \eqref{extra2} that
\begin{equation}\label{eq-78}
|y(s,t)-g_j(x(s,t))|\le C\rho+p_j(|x(s,t)|).
\end{equation}
Using \eqref{eq-76} and \eqref{eq-77} we find that
\begin{equation*}
\begin{split}
\pa_ty(s,t)-&(\pa_tx(s,t))g_j'(x(s,t))=\rho-(y(s,0)-x(s,0))
+(\pa_tx(s,t))p_j'(x(s,t))\\
&=\rho-(y(s,0)-x(s,0))+(x(s,1)-x(s,0))p_j'(x(s,t))\\
&=\rho-(y(s,0)-x(s,0))-x(s,0)p_j'(x(s,t))+x(s,1)p_j'(x(s,t))\\
&\ge (1-C\de)\rho+x(s,1)p_j'(x(s,t))\ge (1-2C\de)\rho+x(s,1)p_j'(|x(s,t)|)\\
&\ge (1-2C\de)\rho+C^{-1}p_j(|x(s,t)|).
\end{split}
\end{equation*}
We now choose $\de$ small
enough so that e.g.~$2C\de<1/2$. This gives the inequality
\begin{equation}\label{eq-79}
\pa_ty(s,t)-(\pa_tx(s,t))g_j'(x(s,t))\ge C^{-1}(\rho+p_j(|x(s,t)|)).
\end{equation}
Combining \eqref{eq-79} with \eqref{eq-73}, \eqref{eq-74}, \eqref{eq-75} 
and \eqref{eq-78} we get
\begin{equation*}
\begin{split}
\pa_th_j(\ga(s,t))&\ge \be_j(\ga(s,t))(\pa_ty(s,t)-(\pa_tx(s,t))g_j'(x(s,t)))\\
&\quad\quad
-|(y(s,t)-g_j(x(s,t)))\pa_t\be_j(\ga(s,t))|\\
&\ge C^{-2}(\rho+p_j(|x(s,t)|))-C^2a(\rho+p_j(|x(s,t)|))\\
&=(C^{-2}-C^2a)(\rho+p_j(|x(s,t)|)).
\end{split}
\end{equation*}
Taking $a<C^{-4}/2$ we obtain a positive bound from below for the right-hand 
side in the last expression, which completes the proof of the lemma.
\end{proof}

\section{Examples and Further Problems}\label{s-6}
         
\subsection{The Necessity of Non-degeneracy Assumptions}\label{s-61}

If one of the cones $\sigma_i(p)$ in~\eqref{eq-cone} is a line it may happen 
that $W(p)\setminus\{p\}$ 
is the union of two components $W(p)_{l}$ and $W(p)_{r}$, each 
bounded by level curves as 
above. In this case there might be several different subharmonic $PH$ 
functions that satisfy condition (i) in Theorem~\ref{main-r1}, as shown by 
Example \ref{ex1} below. Hence something like condition (ii) is 
indeed necessary in order to obtain the conclusion of the aforementioned 
theorem.

\begin{Ex}\label{ex1} 
Set $H_{1}(x,y)=0$, $H_{2}(x,y)= 4x+x^2-y^2$, and $H_{3}(x,y)=-x$. There are 
three level curves through $(0,0)$ to functions of the form $H_i-H_j$ with 
$i\neq j$. These are depicted in Figure 1. Let 
$\varphi=\max\{ H_1\equiv 0,H_2,H_3\}$.
The functions in the figure closest to the origin in each sector are the 
restriction of $\varphi$ to that sector.

\begin{figure}[!htb]\label{fig1}
\centerline{\hbox{\epsfysize=5cm\epsfbox{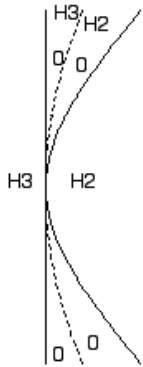}}}
\caption{A non-maximal subharmonic $PH$ function.}
\end{figure}

If one instead defines $\Psi(x,y)$ by changing the value in the two upper 
sectors from $0$ to $H_3 $ respectively $H_2$ then one obtains
a different continuous $PH$ function that is again subharmonic. 
Clearly, every 
neighborhood of the origin still has the property that $\Psi$ is 
equal to each of the three harmonic functions in some subset of positive 
Lebesgue measure. 
%(This simply 
%follows from the fact that subharmonicity is a local property.) 
So $\Psi$ is a maximum of harmonic functions along the curves, hence trivially subharmonic 
away from the origin.  Letting $0\leq \chi\in C_0^\infty(\bR) $ be equal to $1$ near the origin and $\chi_\epsilon(z):=\chi(z/\epsilon),$ this implies that $(1-\chi_\epsilon)\Delta \Psi\geq  0$ in $\D'$. But clearly  $\chi_\epsilon\Delta \Psi\to 0$ in $\D'$ as $\epsilon\to 0$ since
$\Psi=O(\vert z\vert)$. Hence $\Psi $ is subharmonic.
%g.~by 
%calculating integrals along small circles around the origin and checking that 
%these integrals are positive.
\end{Ex}
   
\subsection{On Global Descriptions}\label{s-62}   

In this paper we have only considered the problem of locally characterizing 
the maximum of a finite number of harmonic functions. A natural question is 
to study various situations when a subharmonic $PH$ function is globally the 
maximum of a finite number of harmonic functions. Such a situation occurs 
for instance in \cite{BR}, where the given harmonic functions are linear.
The same conclusion holds when the number of given harmonic functions 
is two as well as in certain other cases. We discuss some of these cases 
in the following examples, which were inspired by~\cite{Mel}.
           
\begin{Ex}
Let $A_1$ and $A_2$ be entire functions such that $A_1(z)\neq A_2(z)$, 
$z\in\bC$
and assume that 
$\Phi:=\chi_1 A_1+\chi_2A_2$ satisfies $\partial\Phi/\partial\bar z\geq 0$, 
where $\chi_1$ and $\chi_2$ are the characteristic functions of the sets 
$M_1$ and $M_2$, respectively (cf.~Notation~\ref{not-r1}).  The first assumption implies that $H_i(z)=\Re\left[\int _0^zA_i(w)dw\right]$, $i=1,2$ are well-defined functions in $\bC$ and that there are no singular points for $H_1-H_2$. For simplicity assume further that level curves to $H_1-H_2$ as well as  the support $\partial\Phi/\partial\bar z$ are connected.
%We further assume that each neighborhood of $p\in \bC$ intersects both 
%$M_1$ and $M_2$ in a set of positive Lebesgue measure. 
If $p\in\overline{M}_1\cap \overline{M}_2$, it follows from 
Theorem~\ref{main-r1} (condition (ii) there being vacuous 
in this case) that there exists a neighborhood $N$ of $p$ and 
constants $c_1(p)$, $c_2(p)$ such that 
$$\Phi=2\frac{\partial}{\partial z}\max\left(H_1+c_1(p),H_2+c_2(p)\right)=2\frac{\partial}{\partial z}\max\left(H_1,H_2+c_2(p)-c_1(p)\right)$$
 In particular, 
the common boundary of $M_1$ and $M_2$ in $N$ is the level curve 
 $H_1-H_2=c_2(p)-c_1(p)$ and this is also the support of $\partial\Phi/\partial\bar z$ in $N$. 
 %Assume for simplicity that the support of 
 %$\partial\Phi/\partial\bar z$ in $\bC$ is connected. 
The local information implies, by the connectedness assumptions,  that globally $c_2(p)-c_1(p)$ is a constant $c$ independent of $p$, and that the support actually consists of the level curve $H_1-H_2=c$, and finally that
$$\Phi=2\frac{\partial}{\partial z}\max\left(H_1,H_2+c\right).$$ 
\end{Ex}

\begin{Ex}
This example is essentially one-dimensional. Assume that
$$\bR=\bigcup_{j=1}^{r}\bar I_j,$$
where the $I_j$ are open pairwise disjoint
intervals. Set $M_j=I_j\times \bR$, $1\le j\le r$, and let $\chi_j(x)$ be the
characteristic function of $I_j$, which we also view as the characteristic
function of $M_j$. Let $h_j(x+\sqrt{-1}y)=a_jx+b_j$, $1\le j\le r$, be
linear functions
on $\bC$ and assume as usual that
$$\chi:=\frac{\partial}{\partial\bar z}\left[\sum_{j=1}^r\frac{\partial
h_j(z)}{\partial z}\chi_j\right]=\sum_{j=1}^r\frac{\partial h_j(z)}{\partial z}
\frac{\partial \chi_j}{\partial\bar  z}=\sum_{j=1}^r \frac{a_j}{2}
\frac{\partial \chi_j}{\partial \bar z}\geq 0.$$
Since $\frac{\partial \chi_j}{\partial \bar z}=\frac{1}{2}\frac{\partial
\chi_j}{\partial x}$ 
we deduce that $\sum_{j=1}^{r}a_j\chi_j$ is an increasing
function of $x$ and thus 
$h(x)=\int_0^x\sum_{j=1}^{r}a_j\chi_j$
is a convex function. Set
\begin{equation*}
H(x,y)=h(y)+h'(y+0)(x-y).
\end{equation*}
By convexity we have
\begin{equation}\label{eq:conv}
h(x)\geq H(x,y),\quad x,y\in \bR,
\end{equation}
with equality when $y=x$. The functions $H(x,y)$ viewed as linear functions
of $x\in \bR$ are independent of $y$ when $y\in I_j$. We denote their common
value for $y\in I_j$ by $\tilde h_j(x)$ and notice that
$\tilde h_j-h_j=C_j$, where $C_j$ is a constant. It follows
from~\eqref{eq:conv} that
$$h(x)=\max_{1\le k\le r}\tilde h_k(x)\text{ in }M_j$$
and then differentiation implies that
$$h'(x)=\frac{\partial}{\partial x}\max_{1\le k\le r}\tilde h_k(x)
=\frac{\partial}{\partial x}\max_{1\le k\le r}\left(h_k(x)+C_k\right).$$
This means precisely that the $PA$ function $\chi$ satisfies
$$\chi=2\frac{\partial}{\partial z}\max_{1\le j\le
r}\left(h_j(z)+C_j\right)$$
and is therefore globally the maximum of a finite number of harmonic
functions.
\end{Ex}

\subsection{Related Questions}\label{s-63}

Let us finally formulate and discuss some interesting related problems. 

\begin{Prob}\label{pb1}
At the moment we do not know although we strongly suspect that locally there 
are in fact only a finite number of possibilities for $\Psi$ even when 
conditions (i)--(iii) are weakened in Theorem \ref{main-r1}. This 
holds e.g.~for the 
function constructed in Example~\ref{ex1}. In 
particular, it seems likely that there always exists a sufficiently 
small neighborhood of $p$ that can be dissected into sectors bounded by 
level curves to $H_{i}-H_{j}$ such that $\Psi$ is constant in each such 
sector. Example~\ref{ex1} suggests that the local behavior of a $PH$ 
subharmonic function is determined by the geometry of the level 
curves $\Gamma_{i,j,k}$ whose study is essentially a problem of a 
combinatorial and topological nature. It would be interesting to give a 
description of this local behavior in terms of Morse theory (the study 
of level curves was Morse's original motivation for his theory, see \cite{K}). 
\end{Prob}

\begin{Prob}\label{pb2}
Another problem is to understand the global behavior of a $PH$ 
subharmonic function and in particular to give criteria saying precisely when 
$\frac{\partial\Psi}{\partial z}$ is the 
derivative of the maximum of a finite number of harmonic functions as in the 
last two examples. This would have interesting applications to uniqueness 
theorems for Cauchy transforms that are algebraic functions as in 
\cite{BR,BBS}.
\end{Prob}

\begin{Prob}\label{pb3}
There are also several connections between the questions
studied in the present paper and the theory of asymptotic solutions to 
differential equations. For instance, sets like those that occur as 
the support of the measures in Theorem~\ref{coro6} play a remarkable role in 
the latter theory (\cite{Fe,K,Wa2,Wa1,Sib}). Moreover, many similar 
techniques are used, e.g.~the admissible sets in \cite{Fe,K} are closely 
related to (though not exactly the same as) the sets $V(z)$ in Lemma 3 above. 
These connections are quite close in the cases studied in \cite{BR,BBS} 
(as well as other cases) and certainly deserve further investigation in view 
of their important applications.
\end{Prob}

\begin{Prob}\label{pb4}
Let $U$ be a domain in $C^n$, where $n\ge 1$. By analogy with 
Definition~\ref{d-r1} and Notation~\ref{not-r1} one can define 
the notions of $PH_n$ and $PA_n$ functions in $U$ as natural 
higher-dimensional generalizations of the concepts of $PH$ and $PA$ functions, 
respectively. It seems reasonable to conjecture that appropriate 
higher-dimensional analogue of Theorem~\ref{main-r1}  
hold for the class $PA_n$ and that
as a consequence one would get a natural extension of 
e.g.~Corollary~\ref{loc-char} to the class $PH_n$.
\end{Prob}

\section*{Appendix. Comments on Some Properties and Definitions}\label{s-ap}

As before, $\chi_\Om$ denotes 
the characteristic function of a set $\Om\subset \bC$ (or $\bR^2$). 
Let us introduce the following additional condition:
an open set $\Om\subset \bR^2$ is said to have {\em property} (*) 
if $\pa\Om$ is of Lebesgue measure $0$ and 
$\pa_z\chi_{\Om},\pa_y\chi_{\Om}$ are measures.

\begin{Lemma}\label{ar-l1}
If $\Om_1,\Om_2\subset \bR^2$ 
have property {\em (*)} then so does $\Om_1\cap\Om_2$.
\end{Lemma}

\begin{proof}
It is clear that $\pa(\Om_1\cap\Om_2)$ is Lebesgue negligible. Let 
$K\subset\bR^2$ be any compact set, choose $\eta\in C_0^{\infty}(\bR^2)$ 
with $\int\!\!\int \eta(x,y)dxdy=1$, define 
$\eta_\eps=\eps^{-2}\eta(x/\eps,y/\eps)$ for $\eps\in(0,1)$ and set 
$\chi_{j,\eps}=\chi_j*\eta_\eps$, where 
$\chi_j=\chi_{_{\Om_j}}$, $j=1,2$. 
Then $0\le \chi_{j,\eps}\le 1$, $\chi_{j,\eps}\to\chi_j$ a.e.~as $\eps\to 0$ 
and 
$
||\pa_x\chi_{j,\eps}||_{L^1(K)}=||\eta_\eps*\pa_x\chi_j||_{L^1(K)}\le C_K,
$
where $C_K$ is independent of $\eps$. Since
$
\pa_x(\chi_{1,\eps}\chi_{2,\eps})=\chi_{1,\eps}\pa_x\chi_{2,\eps}+
\chi_{2,\eps}\pa_x\chi_{1,\eps}
$
it follows that if $\phi\in C_0^{\infty}(\bR^2)$ then
\begin{multline*}
\left|\int\!\!\int \chi_{1,\eps}(x,y)\chi_{2,\eps}(x,y)
\pa_x\phi(x,y)dxdy\right|\le \\
\int\!\!\int |\phi(x,y)|(|\pa_x\chi_{1,\eps}(x,y)|
+|\pa_x\chi_{2,\eps}(x,y)|)
dxdy\le 2C_K||\phi||_{L^{\infty}}.
\end{multline*}
When $\eps\to 0$ this shows that
$$
\left|\int\!\!\int \chi_{1}(x,y)\chi_{2}(x,y)\pa_x\phi(x,y)dxdy\right|\le 
2C_K||\phi||_{L^{\infty}}
$$
and thus $\pa_x(\chi_1\chi_2)$ is a distribution of order $0$ (which extends 
to a measure). This finishes the proof since $\pa_y(\chi_1\chi_2)$ can be 
dealt with in the same way.
\end{proof}

Lemma \ref{ar-l1} shows that 
if we define sets $P^{*}X$ of functions ``piecewise* in $X$'' as in 
Definition \ref{d-r1} 
by demanding in addition that all sets $M_i$ 
have property (*) then $P^{*}X$ are again vector spaces.

\begin{Lemma}\label{ar-l2}
If $u\in P^*X$ is continuous then $\pa_x u,\pa_y u\in P^*X$, where derivatives 
are taken in the distribution sense.
\end{Lemma}

\begin{proof}
Let us write $u=\sum_{i=1}^{r}u_i\chi_i$, where $\chi_i$ is the characteristic 
function of the (open) set $M_i$, $\sum_{i=1}^{r}\chi_i=1$ a.e.~and 
$\pa_x\chi_i,\pa_y\chi_i$ are measures, $1\le i\le r$. Since $u$ is continuous 
we can find $u_\eps\in C^{\infty}(U)$ tending uniformly to $u$ on every 
compact set as $\eps\to 0$. Now
\begin{equation*}\tag{A1}
\begin{split}
\pa_x u_\eps&=\sum_{i=1}^r(\pa_x u_i)\chi_i+\sum_{i=1}^{r}
\chi_i\pa_x(u_\eps-u_i)\\
&=\sum_{i=1}^{r} (\pa_x u_i)\chi_i+
\pa_x\!\left(\sum_{i=1}^{r}(u_\eps-u_i)\chi_i\right)
-\sum_{i=1}^{r}(u_\eps-u_i)\pa_x\chi_i.
\end{split}
\end{equation*}
For every $i$ one has $u_\eps-u_i=u_\eps-u$ in a dense subset of $M_i$. 
It follows that 
$u_\eps\to u_i$ uniformly on every compact subset of $\overline{M}_i$, hence 
also on every compact subset of the support of the measure $\pa_x\chi_i$. 
Therefore, $(u_\eps-u_i)\pa_x\chi_i\to 0$ in $\cD'(\bR^2)$ 
as $\eps\to 0$. This is true for $(u_\eps-u_i)\chi_i$ as well 
and so by letting 
$\eps\to 0$ in (A1) we conclude that 
$\pa_x u_\eps=\sum_{i=1}^{r} (\pa_x u_i)\chi_i$. The same argument applies to 
$\pa_y u$.
\end{proof}

Given a domain $U\subset \bC$ let $S(U)$ be the class of subharmonic 
functions in $U$. Recall Notation \ref{not-Sigma} from \S \ref{s-1}, where 
we already noted the 
(well-known) fact that $\pa_z\phi\in\Si(U)$ whenever 
$\phi\in S(U)$. For completeness we give here a proof of a 
(also well-known) partial converse to this statement.

\begin{Lemma}\label{ref-l-bla}
If $U$ is simply connected and $f\in\Si(U)$ then $f=\pa_z\phi$
for some $\phi\in S(U)$ which is uniquely determined 
modulo an additive constant.
\end{Lemma}

\begin{proof}
Since the operator $\pa_z$ is elliptic we may write $f=\pa_z w$, where 
$w=u+iv\in\cD'(U)$ (cf., e.g., \cite{H}). We get
$
\De u+i\De v=\De w=4\pa_{\bar{z}}\pa_z w=4\pa_{\bar{z}}f\ge 0,
$
which implies that $u\in S(U)$, $v\in H(U)$, and thus $f=\pa_z u+g$, where 
$g=i\pa_z v\in A(U)$. Let $G\in A(U)$ be such that $G'(z)=g(z)$ and define 
$\phi=u+G+\bar{G}$. Then $\phi\in S(U)$ and
$
\pa_z\phi=\pa_z u+\pa_z G=\pa_z u+g=f.
$
The last assertion in the lemma 
follows from the fact that a function $h$ in $U$ is 
constant whenever $h=\bar{h}$ and $\pa_z h=0$.
\end{proof}

\section*{Acknowledgements}

We would like to thank Jan-Erik Bj\"ork and Anders Melin 
for stimulating discussions and useful comments. 
We are especially grateful to the anonymous referee for his detailed 
reports (articles in their own right!) with numerous insightful suggestions 
and an alternative approach for deriving results  
similar to Theorem~\ref{main-r1} under some mild extra assumptions 
(see Theorem~\ref{th-ref} in \S \ref{s-ref}). With his kind permission we 
reproduced large parts of his reports in \S \ref{s-1}, \S \ref{s-ref} 
and the Appendix.

\end{document}